\documentclass{amsart}

\usepackage[english]{babel}
\usepackage{fullpage}
\usepackage{hyperref}

\usepackage{latexsym,amsfonts,amsmath,amssymb,amsthm,amscd,stmaryrd}

\usepackage{yhmath}

\usepackage[dvips]{epsfig}

\usepackage{pdfsync}

\usepackage{graphicx,graphics,psfrag}
\usepackage{float}
\usepackage{color}
\usepackage[table]{xcolor}
\usepackage{subfigure}
\usepackage{graphics,psfrag}

\definecolor{fondrouge}{rgb}{1,0.3,0.3}
\definecolor{fondvert}{rgb}{0.3,1,0.5}

\usepackage{tikz}
\usetikzlibrary{patterns}

\definecolor{rouge}{rgb}{1,0.3,0.3}
\definecolor{vert}{rgb}{0,0.7,0.4}

\newcommand{\bblanc}{\vbox to 7pt{\hbox{
\begin{tikzpicture}[scale=0.6]
\draw [dashed] (3,1) rectangle (3.5,1.5);
\end{tikzpicture}
}}}
\renewcommand{\square}{\bblanc}

\newcommand{\bnoir}{\vbox to 7pt{\hbox{
\begin{tikzpicture}[scale=0.6]
\draw [fill=black!40] (5,1) rectangle (5.5,1.5);
\end{tikzpicture}
}}}

\newcommand{\bleft}{\vbox to 7pt{\hbox{
\begin{tikzpicture}[scale=0.6]
\draw [fill=black!40] (4,1) rectangle (4.5,1.5);
\draw [->,very thick] (4,1.25) -- (4.5,1.25);
\end{tikzpicture}
}}}

\newcommand{\bright}{\vbox to 7pt{\hbox{
\begin{tikzpicture}[scale=0.6]
\draw [fill=black!40] (6,1) rectangle (6.5,1.5);
\draw [<-,very thick] (6,1.25) -- (6.5,1.25);
\end{tikzpicture}
}}}

\newcommand{\G}{\mathbb{G}}
\newcommand{\Z}{\mathbb{Z}}
\newcommand{\N}{\mathbb{N}}
\newcommand{\U}{\mathbb{U}}
\renewcommand{\S}{\mathbb{S}}
\newcommand{\supp}{\mathrm{supp}}
\newcommand{\A}{\mathcal{A}}
\newcommand{\B}{\mathcal{B}}
\newcommand{\E}{\mathcal{E}}
\newcommand{\M}{\mathcal{M}}

\newcommand{\T}{\mathbf{T}}

\newcommand{\s}{\sigma}

\newcommand{\Factor}{\mathbf{Fact}}
\newcommand{\FT}{\mathbf{FT}}
\newcommand{\Product}{\mathbf{Prod}}
\newcommand{\SA}{\mathbf{SA}}

\newcommand{\factor}[2]{\mathbf{Fact}_{#1}\left(#2\right)}
\newcommand{\ft}[2]{\mathbf{FT}_{#1}\left(#2\right)}
\newcommand{\product}[1]{\mathbf{Prod}\left(#1\right)}

\newcommand{\sa}[2]{\mathbf{SA}_{#1}\left(#2\right)}

\newcommand{\TM}{\mathcal{M}}
\newcommand{\SFT}{\mathcal{SFT}}
\newcommand{\sofic}{\mathcal{S}ofic}
\newcommand{\shift}{\mathcal{S}}
\newcommand{\reshift}{\mathcal{RE}}
\newcommand{\fshift}{\mathcal{FS}}
\newcommand{\Uclass}{\mathcal{U}}

\newcommand{\Lang}{\mathcal{L}}

\makeatletter
\renewcommand\section{\@startsection
	{section}
{3}
{0pt}
{-3.5ex plus -1ex minus -.2ex}
{2.3ex plus.2ex}
{\centering\normalfont\Large\scshape}}

\renewcommand\subsection{\@startsection
	{subsection}
{2}
{0pt}
{-1.5ex plus -1ex minus -.2ex}
{0.8ex plus.2ex}
{\normalfont\large\bfseries}}

\renewcommand\subsubsection{\@startsection
	{subsubsection}
{3}
{0pt}
{-1ex plus -1ex minus -.2ex}
{0.6ex plus .2ex}
{\normalfont\itshape}}

\makeatother

\theoremstyle{plain}
\newtheorem{theorem}{Theorem}[section]
\newtheorem{proposition}[theorem]{Proposition}
\newtheorem{fact}[theorem]{Fact}

\theoremstyle{definition}
\newtheorem*{definition}{Definition}

\theoremstyle{remark}
\newtheorem*{remark}{Remark}
\newtheorem{example}{Example}[section]

\newcounter{claimcount}[theorem]

\newcommand{\bprf}[1][Proof:]{\begin{list}{}    {\setlength{\leftmargin}{0.5em}
\setlength{\rightmargin}{0em}  \setlength{\listparindent}{1em}}   \item {\em
\hspace{-1em}  #1  }}
\newcommand{\eprf}{\end{list}}

\usepackage{marvosym,fancybox}

\title[Simulation of effective subshifts by two-dimensional SFT]{Simulation of effective subshifts \\ by two-dimensional subshifts of finite type }
\date{}
\author{Nathalie Aubrun}
\address{LIP,  ENS Lyon, 46 Allée d'Italie Lyon, 69364 Lyon Cedex, France }
\email{nathalie.aubrun@ens-lyon.fr}
\urladdr{http://perso.ens-lyon.fr/nathalie.aubrun/}

\author{Mathieu Sablik}
\address{LATP, Aix-Marseille Universit\'e, 39, rue F. Joliot Curie,13453 Marseille Cedex 13, France }
\email{sablik@latp.univ-mrs.fr}
\urladdr{http://www.latp.univ-mrs.fr/$\sim$sablik/}


\begin{document}

\maketitle

\begin{abstract}
In this article we study how a subshift can simulate another one, where the
notion of simulation is given by operations on subshifts inspired by the
dynamical systems theory (factor, projective subaction...). There exists a
correspondence between the notion of simulation and the set of forbidden
patterns. The main result of this paper states that any effective subshift of
dimension $d$ -- that is a subshift whose set of forbidden patterns can be
generated by a Turing machine -- can be obtained by applying dynamical
operations on a subshift of finite type of dimension $d+1$ -- a subshift that
can be defined by a finite set of forbidden patterns. This result improves
Hochman's~\cite{hochman2007drp}. 
\end{abstract}

\section*{Introduction}\label{section.introduction}

A subshift of dimension $d$ is a closed and shift-invariant subset of
$\A^{\Z^d}$ where $\A$ is a finite alphabet. A subshift can be characterized by
either its language or by a set of forbidden patterns. With this last point of
view, the simplest class is the set of subshifts of finite type, which are
subshifts that can be characterized by a finite set of forbidden patterns. It is
possible to apply dynamical transformations like factor or projective subaction
on a subshift of dimension $d$, and it seems natural to wonder how they modify
the set of forbidden patterns.

In dimension $1$, the class of subshifts of finite type is well understood. In
particular subshifts of finite type are exactly those whose language is accepted
by a local automaton~\cite{beal1993cs}. Given this result, we are naturally
interested in subshifts with a language given by a finite automaton without the
locality condition. This class is entirely characterized in terms of dynamical
operations: it is the class of sofic subshifts, which can all be obtained as a
factor of a subshift of finite type~\cite{lind1995isd}. Thus each sofic subshift
is obtained by a dynamical transformation of a subshift of finite type.

Multidimensional subshifts of finite type are not well understood. For example,
it is not easy to describe their languages. Moreover, in addition to factors,
there exist other types of dynamical transformations on multidimensional
subshifts: for example a subaction of a $d$-dimensional subshift
consists in taking the restriction of a subshift to a subgroup of $\Z^d$.
Hochman~\cite{hochman2007drp} showed that every $d$-dimensional subshift whose
set of forbidden patterns is recursively enumerable can be obtained by
subaction and factor of a $d+2$-subshift of finite type. The main
result of this article states that any effective subshift of dimension $d$ can
be obtained from a SFT of dimension $d+1$, thanks to a subaction and a
factor operation. This result improves Hochman's~\cite{hochman2007drp} since
our construction decreases the dimension. This problem is referenced
in~\cite{boyle2008} and independently of this work there is solution at this problem in~\cite{durand:fixed-point}.

The idea of the proof in~\cite{hochman2007drp} and in this article is to
construct $\T_{\texttt{Final}}$, a three dimensional subshift of finite type
in~\cite{hochman2007drp} (resp. a two-dimensional  subshift of finite type in
this paper), which realizes a given effective subshift
$\Sigma\subset{\A_{\Sigma}}^{\Z}$ in one direction (assume that $d=1$) after a
projection. Thanks to product operation, $\T_{\texttt{Final}}$ is constituted by
different layers, the first one is constituted by the alphabet $\A_{\Sigma}$ and
can be obtained by a projection $\pi$. Then finite type conditions ensure that
for any $x\in \T_{\texttt{Final}}$, one has $\pi(x)_{\Z\times\{ (i,j)\}}\in
\Sigma$ (resp. $\pi(x)_{\Z\times\{i\}}\in \Sigma$) and all these lines are equal; moreover conditions are not
so restrictive and any configuration of $\Sigma$ can be realized by a
configuration of $\T_{\texttt{Final}}$. We here briefly present  the main ideas
of the proof, so that the reader already has in mind the final goal of technical
constructions presented
in this article. The difficulty is to ensure that no forbidden pattern in
$\Sigma$ appears. Since $\Sigma$ is an effective subshift, its forbidden patterns
can be enumerated by a Turing machine. There are classical techniques to
simulate calculations of a Turing machine thanks to a finite type condition (see
Section~\ref{subsection.sftMt}) and the key point of these techniques is that
calculations are embedded into finite computation zones. Thus, we consider a
Turing machine $\TM_{\texttt{Forbid}}$ which has a double role: it both
enumerates the forbidden patterns of $\Sigma$ and checks that none of these
patterns appear in a particular zone around each computation zone, called the
responsibility zone. However, when a Turing machine is constructed in a
two-dimensional subshift of finite type, as in this article, computation zones
are such that a computation is made on a fractured tape (see
Section~\ref{subsection.substsofic}). Consequently forbidden patterns produced
by $\TM_{\texttt{Forbid}}$ are also written on a fractured tape, and comparing
them with non fractured patterns that appear in ${\A_\Sigma}^{\Z^2}$ is not
trivial. To do so, the machine $\TM_{\texttt{Forbid}}$ calls for a second Turing
machine $\TM_\texttt{Search}$ (see Section~\ref{subsection.mt2}). This machine
$\TM_\texttt{Search}$ is given by $\TM_{\texttt{Forbid}}$ an address located in
its responsibility zone, and answers back the letter of $\A_\Sigma$ that appears
at this address. If a forbidden pattern is detected, the machine
$\TM_{\texttt{Forbid}}$ comes into a special state $q_\text{stop}$, whose
presence is forbidden is the final subshift. This ensures that every row
$x_{\Z\times\{ i\}}$ in the final subshift is a configuration of $\Sigma$. So,
the two final operations one has to apply in order to obtain the subshift
$\Sigma$ consist first in taking the projective subaction on $\Z e_1$, where
$e_1$ is the first vector of the canonical basis of $\Z^2$, and then to erase
any information that do not concern $\Sigma$ -- for instance the construction of
computation zones or the SFTs simulating the behaviour of Turing machines --
thanks to a well-chosen letter-to-letter factor.

The difficulty of this construction presented in this paper is to program Turing machines with different size of computation which exchange information in a
two-dimensional subshift of finite type, similar arguments can be found in~\cite{Myers1974,Hanf,Durand}. We note that the authors of~\cite{durand:fixed-point} prove a
similar result based on Kleene's fixed-point theorem. In
that other proof, they do not  recourse to geometric arguments to
describe the circulation of information between the different levels of
computation. 

The paper is organized as follows: in Section~\ref{section.subshifts} we present
five types of operations (product, factor, finite type, projective subaction and
spatial extension) and we formulate classic results with this formalism. In
Section~\ref{section.computationzones}, we present an important tool to define
runs of a Turing machine with a sofic subshift in dimension 2, which is the
construction of an aperiodic SFT that will contain calculations of a Turing
machine and how to code
communication between those different calculations of a Turing machine. These
tools are used to prove our main result in Section~\ref{section.MainResult}. The
main construction of the proof of Theorem~\ref{theorem.mainresult} is built step
by step and for a better understanding, at the end of each
of these subsections the contribution to the final construction is summed up in
a fact. We do not pretend to give a formal proof for these facts, but we hope it
will clarify our intention.

\section{Subshifts and operations on them}\label{section.subshifts}

In this section we recall some basic definitions on subshifts inspired from
symbolic dynamics. We also present some dynamical operations on subshifts, that
were first introduced by Hochman~\cite{hochman2007drp} and then developed by the
authors in~\cite{AubrunSablik09}. 

\subsection{Tilings and subshifts}~\label{subsection.tilingsubshift}
Let $\A$ be a finite alphabet and $d$ be a positive integer. A
\emph{configuration} $x$ is an element of $\A^{\Z^d}$. Let
$\S$ be a finite subset of $\Z^d$. Denote $x_{\S}$
the \emph{restriction} of $x$ to $\S$. A \emph{pattern} is an element
$p\in\A^{\S}$  and $\S$ is the \emph{support} of $p$, which is
denoted by $\supp(p)$. For all $n\in\N$, we call
$\S_n^d=[-n;n]^{d}$ the \emph{elementary support} of size $n$. A
pattern with support $\S_n^d$ is an \emph{elementary pattern}. We
denote by $\mathcal{E}^d_{\A}=\bigcup_{n\in\N} \A^{[-n;n]^{d}}$ the set
of $d$-dimensional elementary patterns. A \emph{$d$-dimensional language}
$\Lang$ is a subset of $\mathcal{E}^d_{\A}$. A pattern $p$ of support
$\S\subset\Z^d$ \emph{appears} in a configuration $x$ if there
exists $i\in\Z^d$ such that for all $j\in\S$, $p_j=x_{i+j}$,
we denote $p\sqsubset x$.

\begin{definition}
A \emph{co-tile set} is a tuple $\tau=(\A,d,P)$ where $P$ is a subset of
$\mathcal{E}^d_{\A}$ called the \emph{set of forbidden patterns}.

A \emph{generalized tiling} by $\tau$ is a configuration $x$ such that for all
$p\in P$, $p$ does not appear in $x$. We denote by $\T_{\tau}$ the set of
generalized tilings by $\tau$. If there is no ambiguity on the alphabet, we just
denote it by $\T_{P}$.
\end{definition}

\begin{remark}
If $P$ is finite, it is equivalent to define a generalized tiling by allowed
patterns or forbidden patterns, the latter being the usual definition of tiling.
\end{remark}

One can define a topology on $\A^{\Z^d}$ by endowing $\A$ with the
discrete topology, and considering the product topology on $\A^{\Z^d}$.
For this topology,  $\A^{\Z^d}$ is a compact metric space on which
$\Z^d$ acts by translation via $\s$ defined by: 
$$\begin{array}{ccccc}
\s_{\A}^i:&\A^{\Z^d} & \longrightarrow &\A^{\Z^d}&\\
&x&\longmapsto & \s_{\A}^i(x)& \textrm{ such that } \s_{\A}^i(x)_u=x_{i+u} \
\forall u\in\Z^d.
\end{array}$$
for all $i$ in $\Z^d$. This action is called the \emph{shift}.

\begin{definition}~\label{definition.language}
A $d$-dimensional subshift on the alphabet $\A$ is a closed and $\s$-invariant
subset of $\A^{\Z^d}$. We denote by $\shift$ (resp. $\shift_d$,
$\shift_{\leq d}$) the set of all subshifts (resp. $d$-dimensional subshifts,
$d'$-dimensional subshifts with $d'\leq d$). 

Let $\T\subseteq \A^{\Z^d}$ be a subshift. Denote
$\Lang_n(\T)\subseteq\A^{[-n;n]^d}$ the set of elementary patterns of size $n$
which appear in some element of $\T$, and $\Lang(\T)=\bigcup_{n\in\N}
\Lang_n(\T)$ the \emph{language} of $\T$ which is the set of elementary patterns
which appear in some element of  $\T$.
\end{definition}

It is also usual to study a subshift as a dynamical
system~\cite{lind1995isd,kitchens1998sd}, the next proposition shows the link
between the two notions.

\begin{proposition}~\label{definition.tiling-subshift}
The set $\T\subset\A^{\Z^d}$ is a subshift if and only if $\T=\T_{\Lang(\T)^c}$
where $\Lang(\T)^c$ is the complement of $\Lang(\T)$ in $\mathcal{E}_{\A}^d$.
\end{proposition}


A set of patterns $P\subseteq\E^d_{\A}$ is \emph{recursively enumerable} if
there exists an effective procedure for listing the patterns of $P$ (see for
instance~\cite{rogersjr1987trf}).

\begin{definition}
It is possible to define different classes of subshifts according to the set of
forbidden patterns: 
\begin{itemize}
\item For a finite alphabet $\A$ and a dimension $d\in\N$, the subshift
$\T_{(\A,d,\emptyset)}=\A^{\Z^d}$ is the \emph{full-shift} of dimension $d$
associated to $\A$. Denote $\fshift$ the set of all full-shifts (for every
finite alphabet $\A$ and dimension $d$).

\item For a finite alphabet $\A$, a dimension $d\in\N$ and a finite set
$P\subseteq\E^d_{\A}$, the subshift $\T_{(\A,d,P)}$ is a \emph{subshift of
finite type}. Denote $\SFT$ the set of all subshifts of finite type. Subshifts
of finite type correspond to the usual notion of tiling.

\item  For a finite alphabet $\A$, a dimension $d\in\N$ and a recursively
enumerable set $P\subseteq\E^d_{\A}$, the subshift $\T_P$ is an \emph{effective
subshift}. Denote $\reshift$ the set of all effective subshifts.
\end{itemize}
\end{definition}

\subsection{Operations on subshifts}\label{subsection.operation}

In this section we describe five operations on subshifts and use them to define
a notion of simulation of a subshift by another one. Operations are gathered in
two groups depending on which part -- the alphabet $\A$ or the group $\Z^d$ --
of a subshift $\T\subseteq \A^{\Z^d}$ they modify.

	\subsubsection{Simulation of a subshift by another one}~\label{subsubsection.simulation}
An \emph{operation} $op$ on subshifts transforms a subshift or a $n$-tuple of
subshifts into another one; it is a function $op:\shift\to\shift$ or
$op:\shift\times\dots\times\shift\to\shift$ that can depend on a parameter. An
operation is not necessarily defined for all subshifts. We
remark that a subshift $\T$ (resp. a pair of subshifts $(\T',\T'')$) and its
image by an operation $op(\T)$ (resp. $op(\T',\T'')$) do not necessary have the
same alphabet or dimension.

Let $Op$ be a set of operations on subshifts. Let $\Uclass\subset\shift$ be a
set of subshifts. We define the \emph{closure} of $\Uclass$ under a set of
operations $Op$, denoted by $\mathcal{C}l_{Op}(\Uclass)$, as the smallest set
stable by $Op$ which contains $\Uclass$. 

We say that a subshift $\T$ \emph{simulates} a subshift $\T'$ by $Op$ if
$\T'\in\mathcal{C}l_{Op}(\T)$. Thus there exists a finite sequence of operations
chosen among $Op$, that transforms $\T$ into $\T'$. We note it by
$\T'\leq_{Op}\T$. Remark that $\mathcal{C}l_{Op}(\T)= \{ \T'\ :\ 
\T'\leq_{Op}\T\}.$

	\subsubsection{Local
transformations}~\label{subsubsection.localtransformations}
We describe three operations that locally modify a subshift
$\T\subseteq\A^{\Z^d}$. The new subshift resulting from the operation will be a
subset of $\B^{\Z^d}$, where $\B$ is a new alphabet.

\paragraph{Product ($\Product$):}

Let $\T_i\subseteq\A_i^{\Z^d}$ for any $i\in\{1,\dots,n\}$ be $n$ subshifts of
the same dimension, define:
$$\product{\T_1,\dots,\T_n}
=\T_1\times\dots\times\T_n\subseteq(\A_1\times\dots\times\mathcal{A}_n)^{\Z^d}
.$$
One has $\mathcal{C}l_{\Product}(\fshift)=\fshift\textrm{ and }
\mathcal{C}l_{\Product}(\SFT)=\SFT.$

\paragraph{Finite type ($\FT$):}

These operations consist in adding a finite number of forbidden patterns to the
initial subshift. Formally, let $\A$ be an alphabet, $P\subseteq\E_{\A}^d$ be a
finite subset and let $\T\subseteq\A^{\Z^d}$ be a subshift. By
Proposition~\ref{definition.tiling-subshift}, there exists $P'$ such that
$\T=\T_{P'}$. Define:
$$ \ft{P}{\T}= \T_{P\cup P'}.$$

Note that $\ft{P}{\T}$ could be empty if $P$ prohibits too many patterns.
By $\FT$, one lists all operations on subshifts which are obtained by this type
of transformation.

By definition of subshift of finite type, one has
$\mathcal{C}l_{\FT}(\fshift)=\SFT$ and $\mathcal{C}l_{\FT}(\fshift)=\SFT$.

\paragraph{Factor ($\Factor$):}

These operations allow to change the alphabet of a subshift by local
modifications. Let $\A$ and $\B$ be two finite alphabets. A \emph{morphism}
$\pi:\A^{\Z^d}\to\B^{\Z^d}$ is a continuous function which commutes with the
shift action (i.e. $\s^i\circ\pi=\pi\circ\s^i$ for all $i\in\Z^d$). In fact,
such a function can be defined locally~\cite{hedlund1969eaa}: that is to say,
there exists $\U\subset\Z^d$ finite, called \emph{neighborhood}, and
$\overline{\pi}:\A^{\U}\to\B$, called \emph{local function}, such that
$\pi(x)_i=\overline{\pi}(\s^i(x)_{\U})$ for all $i\in\Z^d$. 

Let $\pi:\A^{\Z^d}\to\B^{\Z^d}$ be a factor and $\T\subset\A^{\Z^d}$ be a
subshift, define:
$$\factor{\pi}{\T}=\pi(\T).$$

By $\Factor$, one lists all operations on subshifts which are obtained by this
type of transformation.

Example~\ref{Cex1} shows that $\mathcal{C}l_{\Factor}(\SFT)\ne\SFT$.

\begin{example}[$\mathcal{C}l_{\Factor}(\SFT)\ne\SFT$]\label{Cex1}
Consider the alphabet $\{ 0,1,2 \}^{\Z}$ and define $\T=\T_{\{
00,11,02,21\}}$. The factor $\pi$ such that $\pi(0)=\pi(1)=0$ and $\pi(2)=2$
transforms $\T$ into a subshift: 
$$\pi(\T)=\{ x\in\{ 0,2 \}^{\Z} \ : \ \text{ finite blocks of
consecutive 0  are of even length }\}$$ which is called the \emph{even shift}.
It is known that the even shift is not a subshift of finite type (see Example
2.1.9 of~\cite{lind1995isd}), since one need to exclude arbitrarily large blocks
of consecutive $0$'s of odd lengths to describe it.
\end{example}

\begin{definition}
A sofic subshift is a factor of a subshift of finite type. Thus, the set of
sofic subshifts is $\sofic=\mathcal{C}l_{\Factor}(\SFT)$.
\end{definition}

In \cite{lind1995isd}, it is shown that sofic subshifts of dimension $1$ are
subshift which can be defined with a language of forbidden patterns which is
regular. The characterisation is unknown for multidimensional sofic subshifts.

	\subsubsection{Transformations of the group of the
action}~\label{subsubsection.transformationsgroup}
We describe an operation that modify the group on which the subshift is
defined, thus we change the dimension of the subshift.

\paragraph{Projective Subaction ($\SA$):}

These operations allow to take the restriction of a subshift of $\A^{\Z^d}$
according to a subgroup of $\Z^d$. Let $\G$ be a sub-group of $\Z^d$ freely
generated by $u_1,u_2,\dots,u_{d'}$ ($d'\leq d)$. Let $\T\subseteq\A^{\Z^d}$ be
a subshift, define:
$$\sa{\G}{\T}=\left\{y\in\A^{\Z^{d'}}\ :\  \exists x\in\T \textrm{ such that }
\forall i_1,\dots,i_{d'}\in\Z^{d'},
y_{i_1,\dots,i_{d'}}=x_{i_1u_1+\dots+i_{d'}u_{d'}}\right\}.$$

It is easy to prove that $\sa{\G}{\T}$ is a subshift of $\A^{\Z^{d'}}$. One
denotes by $\SA$ the set of all operations on subshifts which are
obtained by this type of operation.

One verifies that $\mathcal{C}l_{\SA}(\SFT)\ne\SFT$ and
$\mathcal{C}l_{\SA}(\SFT)\ne\sofic$ (see respectively Example~\ref{Cex2} and
Example~\ref{Cex3}).

\begin{example}[$\mathcal{C}l_{\SA}(\SFT)\ne\SFT$]\label{Cex2}
We construct a subshift of finite type $\T\subset {\{0,1,2 \} }^{\Z^2}$
such that the projective subaction of $\T$ on the sub-group $\Delta=\{
(x,y)\in\Z^2 : y=x \}\subseteq \Z^2$ is not of finite type. In this example we
want the subshift that appears on $\Delta$ to be $$\left\{ x\in\{ 0,1,2
\}^{\Z} \ : \ \text{ finite blocks of consecutive $0$'s are of even length}\right\}.$$ Define $\overline{F}$ the following set of allowed patterns of size $4$
($.$ symbol may be 1 or 2 but not 0, blank symbol may be 0,1 or 2):

$$
\begin{array}{cccccc}

\begin{array}{|c|c|c|c|}
\hline
 &  & 2 & 0\\
\hline
 & 1 & 0 & 1\\
\hline
2 & 0 & 2 & \\
\hline
0 & 1 &  & \\
\hline
\end{array}
&
;
&
\begin{array}{|c|c|c|c|}
\hline
 &  & . & .\\
\hline
 & 2 & 0 & .\\
\hline
1 & 0 & 1 & \\
\hline
0 & 2 &  & \\
\hline
\end{array}
&
;
&
\begin{array}{|c|c|c|c|}
\hline
 &  & 1 & 0\\
\hline
 & 2 & 0 & 2\\
\hline
. & 0 & 1 & \\
\hline
. & . &  & \\
\hline
\end{array}
&
;

\end{array}
$$

$$
\begin{array}{ccccccc}

\begin{array}{|c|c|c|c|}
\hline
 &  & . & \\
\hline
 & . & . & .\\
\hline
2 & 0 & . & \\
\hline
0 & 1 &  & \\
\hline
\end{array}
&
;
&
\begin{array}{|c|c|c|c|}
\hline
 &  & . & \\
\hline
 & . &  & .\\
\hline
. & . & . & \\
\hline
0 & . &  & \\
\hline
\end{array}
&
;
&
\begin{array}{|c|c|c|c|}
\hline
 &  & 2 & 0\\
\hline
 & . & 0 & 1\\
\hline
. & . & . & \\
\hline
 & . &  & \\
\hline
\end{array}
&
;
&
\begin{array}{|c|c|c|c|}
\hline
 &  & . & 0\\
\hline
 & . & . & .\\
\hline
. &  & . & \\
\hline
 & . &  & \\
\hline
\end{array}

\end{array}
$$
The alternation of 1 and 2 over and under the diagonal of 0 enables us to control
the parity of 0 blocks. Define $F$ as the set of elementary patterns of size 4
that are not in $\overline{F}$. Then if we denote $\T=\T_F$:
$$\sa{\Delta}{\T}=\left\{ x\in\{ 0,1,2 \}^{\Z}\  : \ \text{  blocks of
consecutive $0$'s are of even length}\right\}$$
which is not a subshift of finite type as explained in Example~\ref{Cex1}.
\end{example}

\begin{example}[$\mathcal{C}l_{\SA}(\SFT)\ne\sofic$]\label{Cex3}

The non finite type subshift constructed in Example~\ref{Cex2} is sofic, but it
is possible to obtain non sofic subshifts. We construct a subshift of finite type
$\T$ such that the projection $\sa{\Delta}{\T}$ on the straight line $y=x$ is
not sofic. It is well known that in dimension 1, sofic subshift are exactly
subshifts whose language --- see Definition~\ref{definition.language} --- is a
regular language~\cite{lind1995isd}. The language $\{ a^n b^n :
n\in\N\}$ is non-regular and so we construct a subshift of finite type
$\T\subseteq \A^{\Z^2}$ and a morphism $\pi:\A^{\Z^2}\rightarrow \{
0,a,b\}^{\Z^2}$ such that the only allowed patterns in $\T'=\pi(\T)$ containing
finite blocks of consecutive  $a$'s or $b$'s are those of the form $2n\times 2n$:

$$
\begin{array}{|cccccccc|}
\hline
& & & & & & & 0\\
& & & & & & b &\\
& & & & & \adots & & \\
& & & & b & & & \\
& & & a & & & & \\
& & \adots & & & & & \\
& a & & & & & & \\
0 & & & & & & & \\
\hline
\end{array}
$$

The principle is to construct patterns of even size and to localize the center
of these patterns to distinguish the $a^n$ part from the $b^n$ part.

Denote $\A=\{ *,a,b,0,1,2,3,4\}$. We construct squares formed by any symbols
except the symbol $0$ which forms a background.The symbols $1,2,3$ and $4$ help
to draw the two diagonals of the square and to distinguish in which quadrant we
are. The symbol $*$ only appears on a diagonal of the square, and the other
diagonal contains the $a^n b^n$ part. The presence of the symbol $0$ everywhere
around a finite figure ensures that the two diagonals cross in their middle,
hence the figure pictured is a square. It is possible to describe a finite set
of patterns where the only finite figures on the background formed by $0$'s which
are allowed are even size squares of the form:

$$
\begin{array}{|cccccccc|}
\hline
 0 & 0 & 0 & 0 & 0 & 0 & 0 & 0 \\
 0 & * & 1 & \dots & \dots & 1 & b & 0\\
 0 & 4 & \ddots & 1 & 1 & \adots & 2 & 0\\
 0 & \vdots & 4 & * & b & 2 & \vdots & 0\\
 0 & \vdots & 4 & a & * & 2 & \vdots & 0\\
 0 & 4 & \adots & 3 & 3 & \ddots & 2 & 0\\
 0 & a & 3 & \dots & \dots & 3 & * & 0\\
 0 & 0 & 0 & 0 & 0 & 0 & 0 & 0 \\
\hline
\end{array} \;\;\;(*)
$$
We do not detail the entire set of allowed patterns, but the reader can easily
deduce the missing patterns from those given below:

$$
\text{Squares center: }
\begin{array}{cc}
\begin{array}{|cccc|}
\hline
 * & 1 & 1 & b\\
 4 & * & b & 2\\
 4 & a & * & 2\\
 a & 3 & 3 & *\\
\hline
\end{array}
&
\begin{array}{|cccc|}
\hline
 0 & 0 & 0 & 0\\
 0 & * & b & 0\\
 0 & a & * & 0\\
 0 & 0 & 0 & 0\\
\hline
\end{array}
\end{array}
$$

$$
\text{Squares diagonals: }
\begin{array}{ccccc}

\begin{array}{|ccc|}
\hline
* & 1 & 1\\
4 & * & 1\\
4 & 4 & *\\
\hline
\end{array}
&
\begin{array}{|ccc|}
\hline
1 & 1 & b\\
1 & b & 2\\
b & 2 & 2\\
\hline
\end{array}
&
\begin{array}{|ccc|}
\hline
* & 2 & 2\\
3 & * & 2\\
3 & 3 & *\\
\hline
\end{array}
&
\begin{array}{|ccc|}
\hline
4 & 4 & a\\
4 & a & 3\\
a & 3 & 3\\
\hline
\end{array}
\\

\end{array}
$$

$$
\text{Squares sides: }
\begin{array}{cccccc}

\begin{array}{|ccc|}
\hline
0 & 0 & 0\\
0 & * & 1\\
0 & 4 & *\\
\hline
\end{array}
&
\begin{array}{|ccc|}
\hline
0 & 0 & 0\\
1 & 1 & 1\\
* & 1 & 1\\
\hline
\end{array}
&
\begin{array}{|ccc|}
\hline
0 & 0 & 0\\
1 & 1 & 1\\
1 & 1 & 1\\
\hline
\end{array}
&
\begin{array}{|ccc|}
\hline
0 & 0 & 0\\
1 & 1 & 1\\
1 & 1 & b\\
\hline
\end{array}
&
\begin{array}{|ccc|}
\hline
0 & 0 & 0\\
1 & b & 0\\
b & 2 & 0\\
\hline
\end{array}

\end{array}
$$

$\dots$ and so on for the three other sides.

$$
\text{Uniform domains: }
\begin{array}{cccccc}
\begin{array}{|ccc|}
\hline
0 & 0 & 0\\
0 & 0 & 0\\
0 & 0 & 0\\
\hline
\end{array}
&
\begin{array}{|ccc|}
\hline
1 & 1 & 1\\
1 & 1 & 1\\
1 & 1 & 1\\
\hline
\end{array}
&
\begin{array}{|ccc|}
\hline
2 & 2 & 2\\
2 & 2 & 2\\
2 & 2 & 2\\
\hline
\end{array}
&
\begin{array}{|ccc|}
\hline
3 & 3 & 3\\
3 & 3 & 3\\
3 & 3 & 3\\
\hline
\end{array}
&
\begin{array}{|ccc|}
\hline
4 & 4 & 4\\
4 & 4 & 4\\
4 & 4 & 4\\
\hline
\end{array}
\\
\end{array}
$$

The only configurations one can construct with these allowed patterns are
configurations of $\A^{\Z^2}$ with $0$ everywhere except in some places
where there are arbitrarily large blocks of the form $(*)$, and the
configurations made of the infinite pattern $(*)$. We denote by $\T$ this
subshift of finite type.

Let $\pi$ denote the letter-to-letter morphism defined by $\pi(x)=0$ for $x\in\{
*,1,2,3,4\}$ and $\pi(a)=a$, $\pi(b)=b$. Suppose that $\sa{\Delta}{\T}$ is
sofic. Since $\mathcal{C}l_{\Factor}(\sofic)=\sofic$ then $\pi(\sa{\Delta}{\T})$
would also be sofic, which is absurd since:
$$\pi(\sa{\Delta}{\T})=\T_{\{ba;0a^mb^n0:m\ne n\}}.$$ So this construction
proves that $\mathcal{C}l_{\SA}(\SFT)\ne\sofic$.
\end{example}

The class of SFT is not stable under projective subaction and the class
$\mathcal{C}l_{\SA}(\SFT)$ is studied in~\cite{pavlovschraudner2009}.
Nevertheless a stable class for this operation is known, it is the class of
effective subshifts. This follows from the fact that projective subactions are
special cases of factors of subactions, and by Theorem 3.1 and Proposition 3.3
of ~\cite{hochman2007drp} which establish that symbolic factors and subactions
preserve effectiveness. That is to say $\mathcal{C}l_{\SA}(\reshift)=\reshift$.

With this formalism, the result of M. Hochman~\cite{hochman2007drp} can be
written:
$$\mathcal{C}l_{\Factor,\SA}(\SFT)=\reshift.$$
More precisely, he proves that $\mathcal{C}l_{\Factor,\SA}(\SFT\cap
\shift_{d+2})\cap\shift_{\leq d}=\reshift\cap\shift_{\leq d}$.

In Theorem~\ref{theorem.mainresult}, we show that
$\mathcal{C}l_{\Factor,\SA}(\SFT\cap \shift_{d+1})\cap\shift_{\leq
d}=\reshift\cap\shift_{\leq d}$. Moreover, there are examples of effective
subshifts which are not sofic so
$\mathcal{C}l_{\Factor}(\SFT\cap
\shift_{d})=\sofic\cap\shift_{d}\ne\reshift\cap\shift_{d}$. 

\section{Computation zones for Turing machines}\label{section.computationzones}

In this section we explain how to construct computation zones for a Turing
machine and how to use them to simulate calculations. A Turing machine is a model
of calculation composed by a finite automaton -- the head of calculation -- that
moves on an infinite tape divided into boxes, each box containing a letter that
can be modified by the head. A precise definition of Turing machine will be
given in Subsection~\ref{subsection.defTM}, and it will be explained how to code
the behaviour of the machine thanks to local rules. The main problem
is that this SFT is not enough to code calculations of the machine, since there
is no rule that ensures the calculation is well initialized. So we need to embed
calculations into specific zones. To make sure that the size of these
computation zones is not a constraint and does not prematurely stop a
calculation, we construct arbitrarily large computation zones with a sofic
subshift in Subsection~\ref{subsection.substsofic} and we implement the local
rules of the Turing machine in these zones in
Subsection~\ref{subsection.sftMt}.

\subsection{Local rules to code the behaviour of a Turing
machine}\label{subsection.defTM}

In this article, we consider Turing machines with some restrictions: the
behaviour of the machine will be simulated
only on the empty word (originally the tape only contains blank symbols
$\sharp$). We also assume that the head cannot go to the left of the initial
position. Note that we can impose these restrictions without loss of generality.
First we recall the formal definition of a Turing machine. Remember that a
Turing machine is a model of calculation composed by a finite automaton -- the
head of calculation -- that can be in different states and moves on an infinite
tape divided into boxes, each box containing a letter that can be modified by
the head.

\begin{definition}
Let $\mathcal{M}=(Q,\Gamma,\sharp,q_0,\delta,Q_F)$ be a Turing machine, where:
\begin{itemize}
\item $Q$ is a finite set of states of the head of calculation; $q_0\in Q$ is
the initial state;
\item $\Gamma$ is a finite alphabet;
\item $\sharp\notin\Gamma$ is the blank symbol, with which the tape is initially
filled;
end of any enumerated word;
\item $\delta : Q\times\Gamma\to Q\times\Gamma\times\{\leftarrow,\cdot\,
,\rightarrow\}$ is the transition function. Given the state of the head of
calculation and the letter it can read on the tape --- which thus depends on the
position of the head of calculation on the tape --- the letter on the tape is
replaced or not by another one, the head of calculation moves or not to an
adjacent box and changes or not of state;
\item $F\subset Q_F$ is the set of final states --- when a final state is
reached, the calculation stops.
\end{itemize}
\end{definition}

\begin{example}\label{example.mt}
We consider the Turing machine $\M_\texttt{ex}$ that enumerates on its tape the
words $ab,aabb,aaabbb,\dots$ and never halts. This machine uses the three
letters alphabet $\{a,b,\parallel \}$ and five states $Q=\{ q_0, q_\texttt{a+},
q_\texttt{b+}, q_\texttt{b++},q_\parallel \}$. A separation symbol $\parallel$
is written at the end of each $a^n b^n$. The transition function
$\delta_\texttt{ex}$ is
$$
\begin{array}{cl}
\left.\begin{array}{l}
\delta_\texttt{ex}(q_0,\sharp)=(q_\texttt{b+},a,\rightarrow)\\
\delta_\texttt{ex}(q_\texttt{b+},\sharp)=(q_\parallel,b,\rightarrow)\\
\delta_\texttt{ex}(q_\parallel,\sharp)=(q_\parallel,\parallel,.)
\end{array}\right\}
&
\begin{minipage}{9cm}
Initialization of the tape: the machine writes the first word $ab$ on the tape
and place the head on the separation symbol $\parallel$ to the right of the
word.
\end{minipage}
\\
 & \\
\left.\begin{array}{l}
\delta_\texttt{ex}(q_\parallel,\parallel)=(q_\parallel,\parallel,\leftarrow)\\
\delta_\texttt{ex}(q_\parallel,b)=(q_\parallel,b,\leftarrow)\\
\delta_\texttt{ex}(q_\parallel,a)=(q_\texttt{a+},a,\rightarrow)
\end{array}\right\}
&
\begin{minipage}{9cm}
Suppose some word $a^n b^n\parallel$ is written on the tape, and that the head
is on the $\parallel$ symbol in state $q_\parallel$. The machine looks for the
rightmost symbol a in $a^nb^n$. 
\end{minipage}
\\
& \\
\left.\begin{array}{l}
\delta_\texttt{ex}(q_\texttt{a+},b)=(q_\texttt{b++},a,\rightarrow)\\
\delta_\texttt{ex}(q_\texttt{b++},b)=(q_\texttt{b++},b,\rightarrow)\\
\delta_\texttt{ex}(q_\texttt{b++},\parallel)=(q_\texttt{b+},b,\rightarrow)
\end{array}\right\}
&
\begin{minipage}{9cm}
The machine replaces the leftmost symbol b by a symbol a and looks for the
separation symbol $\parallel$ on the right of the word. Once it has found it, it
is replaced  by $bb\parallel$. The word $a^{n+1}b^{n+1}\parallel$ is now written
on the tape and the head is on the $\parallel$ symbol in state $q_\parallel$.
\end{minipage}
\\
&

\end{array}
$$
A calculation of this machine on an empty tape will always go through the configurations of the tape represented in figure~\ref{figure:comut}.
\begin{figure}[h!]
\begin{center}
\begin{tiny}
$$\begin{array}{c|c|c|c|c|c|c|c|c|c|c|c}

\dots&\,\;\dots\,\;&\,\;\dots\,\;&\,\;\dots\,\;&\,\;\dots\,\;&\,\;\dots\,\;&\,
\;\dots\,\;&\,\;\dots\,\;&\,\;\dots\,\;&\,\;\dots\,\;&\,\;\dots\,\;&\dots\\
\hline
\dots& \sharp  & a  &  a  & a & a & b & \textrm{ $(q_{\texttt{b++}}, b)$}  &
\parallel & \sharp & \sharp & \dots\\
\hline
\dots& \sharp  & a  &  a  & a & a & \textrm{ $(q_{\texttt{b++}}, b)$} & b  &
\parallel & \sharp & \sharp & \dots \\
\hline
\dots& \sharp  & a  &  a  & a &\textrm{ $(q_{\texttt{a+}}, b)$} & b & b  &
\parallel & \sharp & \sharp & \dots \\
\hline
\dots& \sharp  & a  &  a  &\textrm{ $(q_{\parallel},  a)$} & b & b & b  &
\parallel & \sharp & \sharp &\dots \\
\hline
\dots& \sharp  & a  &  a & a  &\textrm{ $(q_{\parallel},  b)$} & b & b  &
\parallel & \sharp & \sharp &  \dots \\
\hline
\dots& \sharp  & a  &  a & a & b  &\textrm{ $(q_{\parallel},  b)$} & b  &
\parallel & \sharp & \sharp &\dots  \\
\hline
\dots& \sharp  & a  &  a & a & b & b  &\textrm{ $(q_{\parallel},  b)$}  &
\parallel & \sharp & \sharp & \dots \\
\hline
\dots& \sharp  & a  &  a & a & b & b  & b &\textrm{ $(q_{\parallel},  \parallel)$} 
& \sharp & \sharp & \dots \\
\hline
\dots& \sharp  & a  &  a & a & b & b  & b &\textrm{ $(q_{\parallel},  \sharp)$} 
& \sharp & \sharp & \dots \\
\hline
\dots& \sharp  & a  &  a & a & b & b & \textrm{ $(q_{\texttt{b+}}, \sharp)$}  &
\sharp & \sharp & \sharp & \dots\\
\hline
\dots& \sharp  & a  &  a & a & b & \textrm{ $(q_{\texttt{b++}}, \parallel)$}  &
\sharp & \sharp & \sharp & \sharp &  \dots \\
\hline
\dots& \sharp  & a  &  a & a &\textrm{ $(q_{\texttt{b++}}, b)$}  &  \parallel &
\sharp & \sharp & \sharp & \sharp & \dots\\
\hline
\dots& \sharp  & a  &  a & \textrm{ $(q_{\texttt{a+}}, b)$} & b  &  \parallel &
\sharp & \sharp & \sharp & \sharp &  \dots \\
\hline
\dots& \sharp  & a  &  \textrm{ $(q_{\parallel}, a)$} & b & b &  \parallel &
\sharp & \sharp & \sharp & \sharp & \dots \\
\hline
\dots& \sharp  & a &a &  \textrm{ $(q_{\parallel}, b)$} & b &  \parallel &
\sharp & \sharp & \sharp & \sharp & \dots\\
\hline
\dots& \sharp  & a &a & b & \textrm{ $(q_{\parallel}, b)$} & \parallel & \sharp
& \sharp & \sharp & \sharp & \dots\\
\hline
\dots& \sharp  & a &a & b& b & \textrm{ $(q_{\parallel}, \parallel)$} & \sharp &
\sharp & \sharp & \sharp &  \dots\\
\hline
\dots& \sharp  & a &a & b& b & \textrm{ $(q_{\parallel}, \sharp)$} & \sharp &
\sharp & \sharp & \sharp &  \dots\\
\hline
\dots& \sharp  & a &a & b& \textrm{ $(q_{\texttt{b+}},\sharp)$}  & \sharp &
\sharp & \sharp & \sharp & \sharp &  \dots\\
\hline
\dots& \sharp  & a &a& \textrm{ $(q_{\texttt{b++}},\parallel)$}  & \sharp &
\sharp & \sharp & \sharp & \sharp & \sharp & \dots \\
\hline
\dots& \sharp  & a & \textrm{ $(q_{\texttt{a+}},b)$} & \parallel  & \sharp &
\sharp & \sharp & \sharp & \sharp & \sharp & \dots \\
\hline
\dots& \sharp  &\textrm{ $(q_{\parallel},a)$} & b & \parallel  & \sharp & \sharp
& \sharp & \sharp & \sharp & \sharp & \dots \\
\hline
\dots& \sharp  &a &\textrm{ $(q_{\parallel},b)$} & \parallel  & \sharp & \sharp
& \sharp & \sharp & \sharp & \sharp & \dots \\
\hline
\dots& \sharp  &a &b &\textrm{ $(q_{\parallel},\parallel)$}  & \sharp & \sharp &
\sharp & \sharp & \sharp & \sharp & \dots \\
\hline
\dots& \sharp  &a &b &\textrm{ $(q_{\parallel},\sharp)$}  & \sharp & \sharp &
\sharp & \sharp & \sharp & \sharp & \dots \\
\hline
\dots& \sharp  &a &\textrm{ $(q_{\texttt{b+}},\sharp)$} & \sharp & \sharp &
\sharp & \sharp & \sharp & \sharp & \sharp & \dots\\
\hline
\dots& \sharp  &\textrm{ $(q_0,\sharp)$} & \sharp & \sharp & \sharp & \sharp &
\sharp & \sharp & \sharp & \sharp & \dots\\
\hline
\end{array}$$
\end{tiny}
\end{center}
\caption{A calculation of this machine on an empty tape will always go through the configurations of the tape.
}\label{figure:comut}
\end{figure}
\end{example}

If an origin is given it is straightforward to describe the behaviour of a
Turing machine with a set of two-dimensional patterns. The first dimension
stands for the tape and second dimension for time evolution. We obtain the {\em
space time diagram of computation of $\TM$} which can be construct locally by
$3\times 2$ allowed patterns: 
\begin{itemize}
\item[$\bullet$] If the pattern codes a part of the tape on which the head of
calculation does not act, the two line of allowed pattern are identical and for
$x,y,z\in\Gamma$ one has:
 $$\begin{array}{|c|c|c|}
\hline
 x & y & z\\
\hline
 x & y & z\\
\hline
\end{array}
$$
 
\item[$\bullet$] If the head of calculation is present in the part of the tape
coded, we code the modification given by the Turing machine. For example the
rule $\delta(q_1,x)=(q_2,y,\leftarrow)$ will be
coded by:$$
\begin{array}{|c|c|c|}
\hline
 (q_2,z) & y & ~z'~\\
\hline
 z & (q_1,x) & ~z'~\\
\hline
\end{array}
$$
\end{itemize}

Denote by $P_{\TM}$ the set of forbidden patterns on the alphabet
$\mathcal{A}_{\TM}=\Gamma\cup(Q\times\Gamma)$ constructed according to the rules
of $\TM$ -- patterns that cannot be seen as coming from the transition function
as above. We can assume that the support of all patterns in $P_{\TM}$ have the
following type: {\vbox to 12pt{\hbox{
\begin{tikzpicture}[scale=0.25]
\draw[dashed] (0,0) rectangle (1,1);
\draw[dashed] (2,0) rectangle (1,1);
\draw[dashed] (2,0) rectangle (3,1);
\draw[dashed] (1,1) rectangle (2,2);
\end{tikzpicture}
}}}. For example, with this assumption the rule
$\delta(q_1,x)=(q_2,y,\leftarrow)$ becomes:
$$
\begin{array}{ccc}
\begin{array}{|c|c|c|}
\hline
 & y &\\
\hline
 ~z~ & (q_1,x) & ~z'~\\
\hline
\end{array}
&
\begin{array}{|c|c|c|}
\hline
  & (q_2,z) &\\
\hline
 ~t~ & z & (q_1,x)\\
\hline
\end{array}
&
\begin{array}{|c|c|c|}
\hline
 & z'& \\
\hline
 (q_1,x) & ~z'~ & ~z''~\\
\hline
\end{array}
\end{array}
$$

Consider now the subshift of finite type $\T_{P_{\TM}}$. It contains an element
that is exactly the space time diagram of computation of $\TM$, but also many
other elements that are inconsistent. With the Turing machine $\TM_\texttt{ex}$
of Example~\ref{example.mt}, the SFT $\T_{P_{\TM_\texttt{ex}}}$ contains an
element where the following configuration of the tape appears
$$
\begin{array}{c|c|c|c|c|c|c|c|c|c|c|c|c|c|c}
\hline
\dots & \sharp & \dots & \sharp & a & b & b & b & b & b &
(q_\parallel,\parallel) & \sharp & \dots & \sharp & \dots\\
\hline
\end{array}
$$
but this configuration is inconsistent since it is never reached by a
calculation of 
$\TM_\texttt{ex}$.

The problem comes from the lack of information about the beginning of a
calculation. We need to specify a point in $\Z^2$ that stands for the origin of
a calculation -- the head of calculation is in the initial state $q_0$, and the
row is filled with blank symbols $\sharp$.

By compactness of the set of configurations of a subshift, it is impossible to
impose that a special symbol appears exactly once in every configuration. 

	\subsection{A substitutive sofic subshift as grid of
computation}\label{subsection.substsofic}

A classical problem in tiling theory is the construction of aperiodic tilings,
that are sets of tiles that can only produce aperiodic configurations. A first
example was initially given by Berger, who proved that the domino problem (is it
possible to tile the whole plane with a given finite set of tiles?) is
undecidable (see \cite{berger1966udp} for the original proof by Berger and
\cite{robinson1971uan} for Robinson's proof with a smaller set of tiles).
Robinson reduces this problem to the Turing machine halting problem, which is
known to be undecidable. The heart of the proof is the construction of an
aperiodic tiling, which codes computation zones for Turing machines. These
computation zones are all finite, but for any calculation of a Turing machine
that stops, it is possible to find a zone large enough that contains it.
Robinson entirely describes a finite set of tiles that produces the tiling, but
there are many techniques to obtain it: Mozes~\cite{mozes1989tss} gives a proof
based on substitutions and Durand, Romashchenko and Shen~\cite{Durand2} propose
a proof based on Kleene's fixed point theorem. We here define computation zones
for Turing machine with a two dimensional substitution.

\paragraph{\texorpdfstring{Definition of the substitution $s_{\texttt{Grid}}$}{Definition of the substitution s}}

Let $\A$ be a finite alphabet. A \emph{$(k,k')$-two dimensional substitution} is
a function $s:\A\rightarrow \A^{\U_{k,k'}}$ where 
$\U_{k,k'}=[0,k-1]\times[0,k'-1]$. We naturally extend $s$ to a function $s^{n,n'} :
\A^{\U_{n,n'}} \rightarrow \A^{\U_{nk,n'k'}}$ by identifying $\A^{\U_{nk,n'k'}}$
with $(\A^{\U_{k,k'}})^{\U_{n,n'}}$.
Starting from a letter placed in $(1,1)\in\Z^2$ and applying
successively $s,s^{k,k'},\dots,s^{k^{n-1},k'^{n-1}}$ we obtain a sequence of
patterns in $\A^{\U_{k^i,k'^i}}$ for $i\in\{0,\dots, n\}$. Such patterns are
called \emph{$s$-patterns}. Note that the substitutions $s$ we define here are
deterministic but one can imagine non deterministic substitutions replacing the
function $s$ by a set of substitution rules, where a letter may have different
images by the substitution. The definition of $s$-patterns naturally extends to
non deterministic substitutions.

\begin{figure}[!ht]\centering
\begin{tikzpicture}[scale=0.8]
\draw (4.8,1.9) node[above]{$\mathcal{G}_1$};
\draw [dashed] (3,1) rectangle (3.5,1.5);
\draw [fill=black!40] (4,1) rectangle (4.5,1.5);
\draw [->,very thick] (4,1.25) -- (4.5,1.25);
\draw [fill=black!40] (5,1) rectangle (5.5,1.5);
\draw [fill=black!40] (6,1) rectangle (6.5,1.5);
\draw [<-,very thick] (6,1.25) -- (6.5,1.25);
\draw (9,2.9) node[above]{$\mathcal{G}_2$};
\draw [color=red,dashed] (7,2) rectangle (7.5,2.5);
\draw [red,line width=2pt] (7.05,2) -- (7.05,2.5);
\draw [red,very thick] (7.25,2) |- (7.5,2.25);
\draw [color=red,dashed] (8,2) rectangle (8.5,2.5);
\draw [red,line width=2pt] (8.05,2) -- (8.05,2.5);
\draw [red,very thick] (8,2.25) -| (8.25,2);
\draw [color=red,dashed] (9,2) rectangle (9.5,2.5);
\draw [red,line width=2pt] (9.05,2) -- (9.05,2.5);
\draw [red,very thick] (9,2.25) -| (9.25,2.5);
\draw [color=red,dashed] (10,2) rectangle (10.5,2.5);
\draw [red,line width=2pt] (10.05,2) -- (10.05,2.5);
\draw [red,very thick] (10.25,2.5) |- (10.5,2.25);
\draw [color=red,dashed] (11,2) rectangle (11.5,2.5);
\draw [red,line width=2pt] (11.05,2) -- (11.05,2.5);
\draw [red,line width=2pt] (11,2.45) -- (11.5,2.45);
\draw [red,very thick] (11.25,2.5) |- (11.5,2.25);
\draw [color=red,dashed] (7,1) rectangle (7.5,1.5);
\draw [red,line width=2pt] (7,1.45) -- (7.5,1.45);
\draw [red,very thick] (7.25,1) |- (7.5,1.25);
\draw [color=red,dashed] (8,1) rectangle (8.5,1.5);
\draw [red,line width=2pt] (8,1.45) -- (8.5,1.45);
\draw [red,very thick] (8,1.25) -| (8.25,1);
\draw [color=red,dashed] (9,1) rectangle (9.5,1.5);
\draw [red,line width=2pt] (9,1.45) -- (9.5,1.45);
\draw [red,very thick] (9,1.25) -| (9.25,1.5);
\draw [color=red,dashed] (10,1) rectangle (10.5,1.5);
\draw [red,line width=2pt] (10,1.45) -- (10.5,1.45);
\draw [red,very thick] (10.25,1.5) |- (10.5,1.25);
\draw [color=red,dashed] (11,1) rectangle (11.5,1.5);
\draw [red,line width=2pt] (11,1.45) -- (11.5,1.45);
\draw [red,line width=2pt] (11.45,1) -- (11.45,1.5);
\draw [red,very thick] (11,1.25) -| (11.25,1.5);
\draw [color=red,dashed] (7,0) rectangle (7.5,0.5);
\draw [red,very thick] (7.25,0) |- (7.5,0.25);
\draw [red,line width=2pt] (7.45,0.5) -- (7.45,0);
\draw [color=red,dashed] (8,0) rectangle (8.5,0.5);
\draw [red,line width=2pt] (8.45,0.5) -- (8.45,0);
\draw [red,very thick] (8,0.25) -| (8.25,0);
\draw [color=red,dashed] (9,0) rectangle (9.5,0.5);
\draw [red,line width=2pt] (9.45,0.5) -- (9.45,0);
\draw [red,very thick] (9,0.25) -| (9.25,0.5);
\draw [color=red,dashed] (10,0) rectangle (10.5,0.5);
\draw [red,line width=2pt] (10.45,0.5) -- (10.45,0);
\draw [red,very thick] (10.25,0.5) |- (10.5,0.25);
\draw [color=red,dashed] (11,0) rectangle (11.5,0.5);
\draw [red,line width=2pt] (11.45,0.5) -- (11.45,0);
\draw [red,line width=2pt] (11,0.05) -- (11.5,0.05);
\draw [red,very thick] (11,0.25) -| (11.25,0);

\draw [color=red,dashed] (7,-1) rectangle (7.5,-0.5);
\draw [red,line width=2pt] (7,-0.95) -- (7.5,-0.95);
\draw [red,very thick] (7.25,-1) |- (7.5,-0.75);
\draw [color=red,dashed] (8,-1) rectangle (8.5,-0.5);
\draw [red,line width=2pt] (8,-0.95) -- (8.5,-0.95);
\draw [red,very thick] (8,-0.75) -| (8.25,-1);
\draw [color=red,dashed] (9,-1) rectangle (9.5,-0.5);
\draw [red,line width=2pt] (9,-0.95) -- (9.5,-0.95);
\draw [red,very thick] (9,-0.75) -| (9.25,-0.5);
\draw [color=red,dashed] (10,-1) rectangle (10.5,-0.5);
\draw [red,line width=2pt] (10,-0.95) -- (10.5,-0.95);
\draw [red,very thick] (10.25,-0.5) |- (10.5,-0.75);
\draw [color=red,dashed] (11,-1) rectangle (11.5,-0.5);
\draw [red,line width=2pt] (11,-0.95) -- (11.5,-0.95);
\draw [red,line width=2pt] (11.05,-0.5) -- (11.05,-1);
\draw [red,very thick] (11.25,-1) |- (11.5,-0.75);
\end{tikzpicture}
\caption{The alphabets $\mathcal{G}_1$ and $\mathcal{G}_2$ on which the
substitution is defined.}
\label{figure.alphabets}
\end{figure}
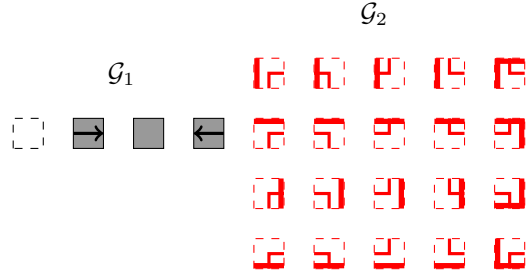

To describe the grid of computation, we consider two alphabets $\mathcal{G}_1$
and $\mathcal{G}_2$ (see Figure~\ref{figure.alphabets}). The alphabet
$\mathcal{G}_1=\{\bblanc,\bnoir,\bleft ,\bright\}$ describes the zones of
computation, $\bnoir$, $\bleft$ and $\bright$ are called {\em computation boxes}
 where the computation holds and $\bblanc$ are called {\em communication boxes}
through which computation boxes can send information. More precisely, $\bleft$
and $\bright$ are called {\em border computation boxes}. The alphabet
$\mathcal{G}_2$ is constituted by lines which describe communication
channels between the different zones of computation.

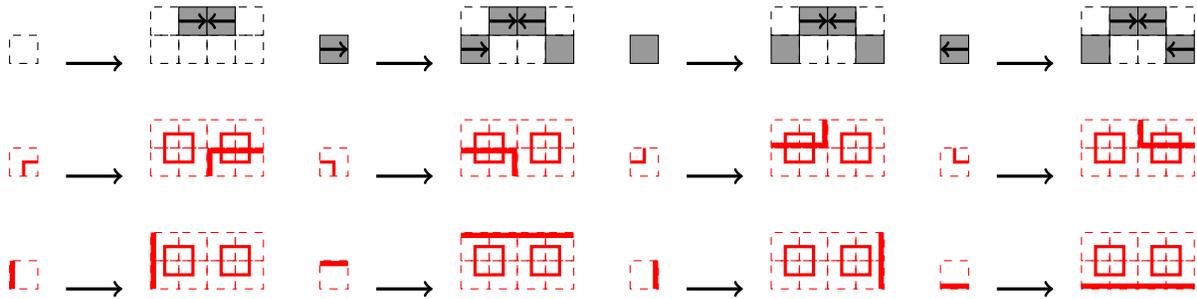
\begin{figure}[!ht]\centering
\begin{tikzpicture}[scale=0.75]
\draw [dashed] (0,9) rectangle (0.5,9.5);
\draw [->,very thick] (1,9) -- (2,9);
\draw [dashed] (2.5,9) rectangle (3,9.5);
\draw [dashed] (2.5,9.5) rectangle (3,10);
\draw [dashed] (3,9) rectangle (3.5,9.5);
\draw [fill=black!40] (3,9.5) rectangle (3.5,10);
\draw [->,very thick] (3,9.75) -- (3.5,9.75);
\draw [dashed] (3.5,9) rectangle (4,9.5);
\draw [fill=black!40] (3.5,9.5) rectangle (4,10);
\draw [<-,very thick] (3.5,9.75) -- (4,9.75);
\draw [dashed] (4,9) rectangle (4.5,9.5);
\draw [dashed] (4,9.5) rectangle (4.5,10);
\draw [fill=black!40] (5.5,9) rectangle (6,9.5);
\draw [->,very thick] (5.5,9.25) -- (6,9.25);
\draw [->,very thick] (6.5,9) -- (7.5,9);
\draw [fill=black!40] (8,9) rectangle (8.5,9.5);
\draw [->,very thick] (8,9.25) -- (8.5,9.25);
\draw [dashed] (8,9.5) rectangle (8.5,10);
\draw [dashed] (8.5,9) rectangle (9,9.5);
\draw [fill=black!40] (8.5,9.5) rectangle (9,10);
\draw [->,very thick] (8.5,9.75) -- (9,9.75);
\draw [dashed] (9,9) rectangle (9.5,9.5);
\draw [fill=black!40] (9,9.5) rectangle (9.5,10);
\draw [<-,very thick] (9,9.75) -- (9.5,9.75);
\draw [fill=black!40] (9.5,9) rectangle (10,9.5);
\draw [dashed] (9.5,9.5) rectangle (10,10);
\draw [fill=black!40] (11,9) rectangle (11.5,9.5);
\draw [->,very thick] (12,9) -- (13,9);
\draw [fill=black!40] (13.5,9) rectangle (14,9.5);
\draw [dashed] (13.5,9.5) rectangle (14,10);
\draw [dashed] (14,9) rectangle (14.5,9.5);
\draw [fill=black!40] (14,9.5) rectangle (14.5,10);
\draw [->,very thick] (14,9.75) -- (14.5,9.75);
\draw [dashed] (14.5,9) rectangle (15,9.5);
\draw [fill=black!40] (14.5,9.5) rectangle (15,10);
\draw [<-,very thick] (14.5,9.75) -- (15,9.75);
\draw [fill=black!40] (15,9) rectangle (15.5,9.5);
\draw [dashed] (15,9.5) rectangle (15.5,10);
\draw [fill=black!40] (16.5,9) rectangle (17,9.5);
\draw [<-,very thick] (16.5,9.25) -- (17,9.25);
\draw [->,very thick] (17.5,9) -- (18.5,9);
\draw [fill=black!40] (19,9) rectangle (19.5,9.5);
\draw [dashed] (19,9.5) rectangle (19.5,10);
\draw [dashed] (19.5,9) rectangle (20,9.5);
\draw [fill=black!40] (19.5,9.5) rectangle (20,10);
\draw [->,very thick] (19.5,9.75) -- (20,9.75);
\draw [dashed] (20,9) rectangle (20.5,9.5);
\draw [fill=black!40] (20,9.5) rectangle (20.5,10);
\draw [<-,very thick] (20,9.75) -- (20.5,9.75);
\draw [fill=black!40] (20.5,9) rectangle (21,9.5);
\draw [<-,very thick] (20.5,9.25) -- (21,9.25);
\draw [dashed] (20.5,9.5) rectangle (21,10);
\draw [dashed,color=red] (0,7) rectangle (0.5,7.5);
\draw [red,very thick] (0.25,7) |- (0.5,7.25);
\draw [->,very thick] (1,7) -- (2,7);
\draw [dashed,color=red] (2.5,7) rectangle (3,7.5);
\draw [dashed,color=red] (2.5,7.5) rectangle (3,8);
\draw [dashed,color=red] (3,7) rectangle (3.5,7.5);
\draw [dashed,color=red] (3,7.5) rectangle (3.5,8);
\draw [dashed,color=red] (3.5,7) rectangle (4,7.5);
\draw [dashed,color=red] (3.5,7.5) rectangle (4,8);
\draw [dashed,color=red] (4,7) rectangle (4.5,7.5);
\draw [dashed,color=red] (4,7.5) rectangle (4.5,8);
\draw [very thick,color=red] (2.75,7.25) rectangle (3.25,7.75);
\draw [very thick,color=red] (3.75,7.25) rectangle (4.25,7.75);
\draw [line width=2pt,red] (3.55,7) |- (4.5,7.45);
\draw [dashed,color=red] (5.5,7) rectangle (6,7.5);
\draw [red,very thick] (5.5,7.25) -| (5.75,7);
\draw [->,very thick] (6.5,7) -- (7.5,7);
\draw [dashed,color=red] (8,7) rectangle (8.5,7.5);
\draw [dashed,color=red] (8,7.5) rectangle (8.5,8);
\draw [dashed,color=red] (8.5,7) rectangle (9,7.5);
\draw [dashed,color=red] (8.5,7.5) rectangle (9,8);
\draw [dashed,color=red] (9,7) rectangle (9.5,7.5);
\draw [dashed,color=red] (9,7.5) rectangle (9.5,8);
\draw [dashed,color=red] (9.5,7) rectangle (10,7.5);
\draw [dashed,color=red] (9.5,7.5) rectangle (10,8);
\draw [very thick,color=red] (8.25,7.25) rectangle (8.75,7.75);
\draw [very thick,color=red] (9.25,7.25) rectangle (9.75,7.75);
\draw [line width=2pt,red] (8.95,7) |- (8,7.45);
\draw [dashed,color=red] (11,7) rectangle (11.5,7.5);
\draw [red,very thick] (11,7.25) -| (11.25,7.5);
\draw [->,very thick] (12,7) -- (13,7);
\draw [dashed,color=red] (13.5,7) rectangle (14,7.5);
\draw [dashed,color=red] (13.5,7.5) rectangle (14,8);
\draw [dashed,color=red] (14,7) rectangle (14.5,7.5);
\draw [dashed,color=red] (14,7.5) rectangle (14.5,8);
\draw [dashed,color=red] (14.5,7) rectangle (15,7.5);
\draw [dashed,color=red] (14.5,7.5) rectangle (15,8);
\draw [dashed,color=red] (15,7) rectangle (15.5,7.5);
\draw [dashed,color=red] (15,7.5) rectangle (15.5,8);
\draw [very thick,color=red] (13.75,7.25) rectangle (14.25,7.75);
\draw [very thick,color=red] (14.75,7.25) rectangle (15.25,7.75);
\draw [line width=2pt,red] (13.5,7.55) -| (14.45,8);
\draw [dashed,color=red] (16.5,7) rectangle (17,7.5);
\draw [red,very thick] (16.75,7.5) |- (17,7.25);
\draw [->,very thick] (17.5,7) -- (18.5,7);
\draw [dashed,color=red] (19,7) rectangle (19.5,7.5);
\draw [dashed,color=red] (19,7.5) rectangle (19.5,8);
\draw [dashed,color=red] (19.5,7) rectangle (20,7.5);
\draw [dashed,color=red] (19.5,7.5) rectangle (20,8);
\draw [dashed,color=red] (20,7) rectangle (20.5,7.5);
\draw [dashed,color=red] (20,7.5) rectangle (20.5,8);
\draw [dashed,color=red] (20.5,7) rectangle (21,7.5);
\draw [dashed,color=red] (20.5,7.5) rectangle (21,8);
\draw [very thick,color=red] (19.25,7.25) rectangle (19.75,7.75);
\draw [very thick,color=red] (20.25,7.25) rectangle (20.75,7.75);
\draw [line width=2pt,red] (20.05,8) |- (21,7.55);
\draw [dashed,color=red] (0,5) rectangle (0.5,5.5);
\draw [->,very thick] (1,5) -- (2,5);
\draw [dashed,color=red] (2.5,5) rectangle (3,5.5);
\draw [dashed,color=red] (2.5,5.5) rectangle (3,6);
\draw [dashed,color=red] (3,5) rectangle (3.5,5.5);
\draw [dashed,color=red] (3,5.5) rectangle (3.5,6);
\draw [dashed,color=red] (3.5,5) rectangle (4,5.5);
\draw [dashed,color=red] (3.5,5.5) rectangle (4,6);
\draw [dashed,color=red] (4,5) rectangle (4.5,5.5);
\draw [dashed,color=red] (4,5.5) rectangle (4.5,6);
\draw [very thick,color=red] (2.75,5.25) rectangle (3.25,5.75);
\draw [very thick,color=red] (3.75,5.25) rectangle (4.25,5.75);
\draw [line width=2pt,red] (0.05,5) -- (0.05,5.5);
\draw [line width=2pt,red] (2.55,5) -- (2.55,6);
\draw [dashed,color=red] (5.5,5) rectangle (6,5.5);
\draw [->,very thick] (6.5,5) -- (7.5,5);
\draw [dashed,color=red] (8,5) rectangle (8.5,5.5);
\draw [dashed,color=red] (8,5.5) rectangle (8.5,6);
\draw [dashed,color=red] (8.5,5) rectangle (9,5.5);
\draw [dashed,color=red] (8.5,5.5) rectangle (9,6);
\draw [dashed,color=red] (9,5) rectangle (9.5,5.5);
\draw [dashed,color=red] (9,5.5) rectangle (9.5,6);
\draw [dashed,color=red] (9.5,5) rectangle (10,5.5);
\draw [dashed,color=red] (9.5,5.5) rectangle (10,6);
\draw [very thick,color=red] (8.25,5.25) rectangle (8.75,5.75);
\draw [very thick,color=red] (9.25,5.25) rectangle (9.75,5.75);
\draw [line width=2pt,red] (5.5,5.45) -- (6,5.45);
\draw [line width=2pt,red] (8,5.95) -- (10,5.95);
\draw [dashed,color=red] (11,5) rectangle (11.5,5.5);
\draw [->,very thick] (12,5) -- (13,5);
\draw [dashed,color=red] (13.5,5) rectangle (14,5.5);
\draw [dashed,color=red] (13.5,5.5) rectangle (14,6);
\draw [dashed,color=red] (14,5) rectangle (14.5,5.5);
\draw [dashed,color=red] (14,5.5) rectangle (14.5,6);
\draw [dashed,color=red] (14.5,5) rectangle (15,5.5);
\draw [dashed,color=red] (14.5,5.5) rectangle (15,6);
\draw [dashed,color=red] (15,5) rectangle (15.5,5.5);
\draw [dashed,color=red] (15,5.5) rectangle (15.5,6);
\draw [very thick,color=red] (13.75,5.25) rectangle (14.25,5.75);
\draw [very thick,color=red] (14.75,5.25) rectangle (15.25,5.75);
\draw [line width=2pt,color=red] (11.45,5.5) -- (11.45,5);
\draw [line width=2pt,color=red] (15.45,6) -- (15.45,5);
\draw [dashed,color=red] (16.5,5) rectangle (17,5.5);
\draw [->,very thick] (17.5,5) -- (18.5,5);
\draw [dashed,color=red] (19,5) rectangle (19.5,5.5);
\draw [dashed,color=red] (19,5.5) rectangle (19.5,6);
\draw [dashed,color=red] (19.5,5) rectangle (20,5.5);
\draw [dashed,color=red] (19.5,5.5) rectangle (20,6);
\draw [dashed,color=red] (20,5) rectangle (20.5,5.5);
\draw [dashed,color=red] (20,5.5) rectangle (20.5,6);
\draw [dashed,color=red] (20.5,5) rectangle (21,5.5);
\draw [dashed,color=red] (20.5,5.5) rectangle (21,6);
\draw [very thick,color=red] (19.25,5.25) rectangle (19.75,5.75);
\draw [very thick,color=red] (20.25,5.25) rectangle (20.75,5.75);
\draw [line width=2pt,color=red] (16.5,5.05) -- (17,5.05);
\draw [line width=2pt,color=red] (19,5.05) -- (21,5.05);
\end{tikzpicture}
\caption{Basic elements to define the substitution rules of $s$. The first row
lists the substitution rules of $s_1$ on the alphabet $\mathcal{G}_1$. The
second and third rows contain substitution rules of $s_2$ on some of the letters
of
$\mathcal{G}_2$. All substitution rules of $s_2$ on $\mathcal{G}_2$ can be
obtained by superimposing a substitution rule of the second row and a
substitution rule of the third row. One can deduce substitution rules of
$s_{\texttt{Grid}}$ on the alphabet $\mathcal{G}_1\times\mathcal{G}_2$ by
superimposing a rule of $s_1$ on $\mathcal{G}_1$ and a rule of $s_2$ on
$\mathcal{G}_2$.}
\label{figure.substitution-rule}
\end{figure}

\begin{figure}[!ht]
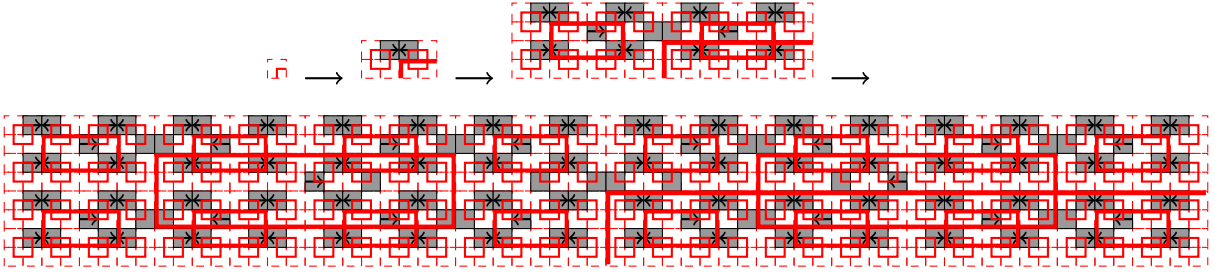
\centering
\include{figure_iteration_bis}
\caption{Four iterations of the substitution $s_{\texttt{Grid}}$ starting from
an element of $\mathcal{G}_1\times\mathcal{G}_2$.}
\label{figure.iteration}
\end{figure}

We define two $(4,2)$-two dimensional substitutions, $s_1$ on $\mathcal{G}_1$
and $s_2$ on $\mathcal{G}_2$ (see Figure~\ref{figure.substitution-rule} for the
substitution rules). Then, we define the product substitution
$s_{\texttt{Grid}}=s_1\times s_2$ on $\mathcal{G}_1\times\mathcal{G}_2$.
Iterations of $s_{\texttt{Grid}}$ on any pattern of
$\mathcal{G}_1\times\mathcal{G}_2$ produce arbitrarily large computation zones
with communication channels between them (this will be detailed in
Section~\ref{subsection.computationzones} and
Section~\ref{subsection.channels}). See Figure~\ref{figure.iteration} for an
example of an iteration of $s_{\texttt{Grid}}$.

We denote by $\pi_{\mathcal{G}_1}$ (resp. $\pi_{\mathcal{G}_2}$) the projection
on $\mathcal{G}_1$ (resp. $\mathcal{G}_2$).

\paragraph{Sofic subshift generated by the substitution}

Given a substitution $s$, recall that a $s$-pattern is a pattern obtained by
iteration of the substitution $s$ on a letter (for instance in
Figure~\ref{figure.iteration} are drawn $s$-patterns obtained after four
iterations on the letter $\square$). The \emph{subshift generated by a
substitution $s$}, denoted $\T_s$, is the set of configurations $x$ such that
any pattern that appears in $x$ also appears in a $s$-pattern.

S. Mozes studied more general substitutions -- non deterministic ones and
substitution rules may be of different sizes -- and proved that if the
substitution $s$ satisfies some good property and has only strictly
two-dimensional substitution rules, then $\T_s$ is a sofic subshift (see Theorem
4.1 of~\cite{mozes1989tss}). In particular Mozes theorem can be applied for all deterministic substitutions, that is to say that all letter have only one image, like $s_{\texttt{Grid}}$. As a consequence, the following
holds

\begin{fact}\label{fact.SoficGrid}
The subshift generated by $s_{\texttt{Grid}}$,
$$\T_{\texttt{Grid}}=\T_{s_{\texttt{Grid}}}=\left\{x\in\left(\mathcal{G}
_1\times\mathcal{G}_2\right)^{\Z^2} :  \textrm{ for all } u\sqsubset x \textrm{
there exists }n\in\N\textrm{ such that }u\sqsubset
s_{\texttt{Grid}}^n({\vbox to
7pt{\hbox{
\begin{tikzpicture}[scale=0.7]
 \draw [dashed] (0,0) rectangle (0.5,0.5);
 \draw [red,very thick] (0.25,0) |- (0.5,.25);
\end{tikzpicture}
}}})\right\}$$
is a two-dimensional sofic subshift.
\end{fact}

\begin{remark}
Note that $\pi_{\mathcal{G}_1}(\T_{\texttt{Grid}})=\T_{s_1}$ and 
$\pi_{\mathcal{G}_2}(\T_{\texttt{Grid}})=\T_{s_2}$ but $\T_{\texttt{Grid}}$ is
different of $\T_{s_1}\times\T_{s_2}$.
\end{remark}

A substitution $s:\A\rightarrow\A^{\U_{k,k'}}$ may be extended into an
application $\widetilde{s}:\A^{\Z^2}\rightarrow\A^{\Z^2}$. This substitution has
\emph{unique derivation} if for every element  $x\in\T_s$ there exists an unique
$y\in\A^{\Z^2}$ and an unique $i\in\U_{k,k'}$ such that $\widetilde{s}(y)=\sigma^i(x)$.

Since the pattern $\begin{tikzpicture}[scale=0.6]
\draw [fill=black!40] (3,9) rectangle (3.5,9.5);
\draw [->,very thick] (3,9.25) -- (3.5,9.25);
\draw [fill=black!40] (3.5,9) rectangle (4,9.5);
\draw [<-,very thick] (3.5,9.25) -- (4,9.25);
\end{tikzpicture}$ appears in each rules of the substitutions $s_1$, for every configuration $x\in\T_{s_1}$, there exists $(i,j)\in[0,3]\times[0,1]$ such that $x_{\{n_1+i+1\}\times[n_2+j+1,n_2+j+2]}= \begin{tikzpicture}[scale=0.6]
\draw [fill=black!40] (3,9) rectangle (3.5,9.5);
\draw [->,very thick] (3,9.25) -- (3.5,9.25);
\draw [fill=black!40] (3.5,9) rectangle (4,9.5);
\draw [<-,very thick] (3.5,9.25) -- (4,9.25);
\end{tikzpicture}$ for all $(n_1,n_2)\in\N\times\N$. Moreover this pattern cannot appear in other position so $(i,j)$ is chosen in an unique way. Consider the plane partition $([n_1+i,n_1+i+3]\times[n_2+j,n_2+j+1])_{(n_1,n_2)\in\N\times\N}$ of the configuration $x$, since all boxes have different image by the substitution, this plane partition gives an unique antecedent y by $\widetilde{s}_1$. We deduce that $\widetilde{s}_1(y)=\s^{i,j}(x)$. Thus $\T_{s_1}$ has unique derivation. The same type of reasoning holds for $\widetilde{s_2}$.

\begin{fact}\label{fact.UniqueDerivation}
The substitutions $s_1$ and $s_2$ have unique derivation. 
\end{fact}


\subsection{Use of communication channels}\label{subsection.communication}

A \emph{communication channel} is a sequence of adjacent boxes marked by a
special
symbol -- we call these marked boxes \emph{channel boxes}. The channel begins
and ends with two computation boxes. In our construction, the channel boxes are
of two types which can appear in the same box:
\begin{itemize}
\item communication boxes $\bblanc$ from alphabet $\mathcal{G}_1$ that will be
used for internal communication and communication between adjacent Turing machines;
\item symbols from alphabet $\mathcal{G}_2$ that will be used for communication
between non adjacent Turing machines.
\end{itemize}

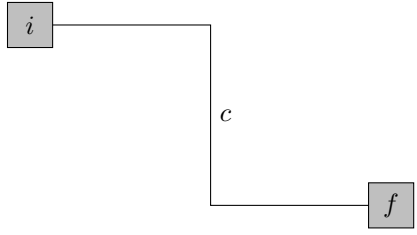
\begin{figure}[!ht]\centering
\begin{tikzpicture}[scale=0.6]
 \draw[fill=gray!50] (0,0)rectangle(1,1);
 \draw (0.5,0.5) node{$i$};
 \draw[fill=gray!50] (8,-4)rectangle(9,-3);
 \draw (8.5,-3.5) node{$f$};
 \draw (1,0.5)--(4.5,0.5)--(4.5,-3.5)--(8,-3.5);
 \draw (4.5,-1.5) node[right] {$c$};
\end{tikzpicture}
\caption{A communication channel denoted $c$ between the computation boxes $i$ and $f$.}
\label{figure.channel}
\end{figure}

A \emph{transfer of information} consists in three objects (see Figure~\ref{figure.channel}):
\begin{itemize}
 \item an \emph{initial computation box} denoted $i$ and a \emph{final
computation box} denoted $f$
 \item local rules that determine the symbol transferred through a channel,
depending on the direction of the channel starting from $i$  (resp.reaching $f$)
and the symbol contained in the box $i$ (resp. in the box $f$)
 \item a communication channel $c$.
\end{itemize}

Two adjacent communication boxes carry the same symbol, which is transferred
through the channel. Note that the computation boxes $i$ and $f$ are not
necessary identical -- for example a rule local may make a change at the end of
the communication channel. The same computation box may be at the extremity of
different communication channels.

Note that a communication box may belong to multiple communication channels, but
this number must be bounded -- in our construction the maximum number of
channels going through a communication box will be 3 -- internal communication
inside a computation zone, communication between two adjacent computation zones
of same level and communication between computation zones of different levels.

\begin{fact}\label{fact.transfer}
Given a subshift $\Sigma$ that contains communication channels, it is possible
to code transfers of information through these channels thanks to a product and
a finite type operations, provided the symbols transferred locally depend on the
symbol contained in the initial and final computation boxes of the channel.
\end{fact}

\subsection{Description of computation zones}\label{subsection.computationzones}

In this section we only consider the $\mathcal{G}_1$ part of the sofic subshift
$\T_{\texttt{Grid}}$. We here describe the grid where computations hold for an element of $\pi_{\mathcal{G}_1}(\T_{\texttt{Grid}})$:
horizontal dimension stands for the tape and vertical dimension for time
evolution. On a horizontal line, a {\em zone of computation} is constituted by a
group of computation boxes $\bnoir$ located between $\bleft $ on the left and
$\bright$ on the right. The {\em size} of a computation zone is the number of
computation boxes which constitute the zone.

\begin{figure}[!ht]
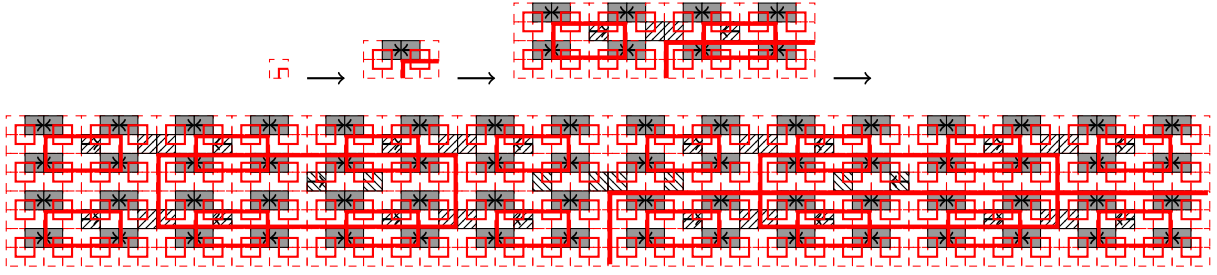
\centering
\include{figure_iteration_ter}
\caption{Four iterations of the substitution $s_{\texttt{Grid}}$. Computation
zones of level 1,2 and 3 appear on the last pattern. The computation boxes of
the first (resp. second and third) level are pictured with plain (resp. hashed
with SW-NE lines and hashed with NW-SE lines) pattern.}
\label{figure.IterationColored}
\end{figure}

Consider a configuration $x\in\T_{\texttt{Grid}}$ and a computation zone in $x$.
Since $s_1$ has unique derivation (see Fact~\ref{fact.UniqueDerivation}) for any
integer $n$, there exists a unique way to partition $x$ into $4^n\times2^n$
rectangles so that each of these rectangles is a $s_1^n(a)$ for some
$a\in\mathcal{G}_1$. So there exists a minimal integer $n$ such that the
computation zone of $x$ appears in $s_1^n(a)$ for some $a\in\mathcal{G}_1$. We
call this integer the \emph{level} of the computation zone (see Figure~\ref{figure.IterationColored}).

At the iteration $n$ of the substitution on $\bblanc$, that is to say
$s^n(\bblanc)$, we  obtained a rectangle of size $4^n\times 2^n$. By
induction, for all $m\in[1,n]$, we  get $4^{n-m}* 2^{n-m}=8^{n-m}$ zones of
computation of level $m$ in $s^n(\bblanc)$, the size of these zones of
computation is $2^m$. More precisely, if on the line $j\in[0,2^n-1]$ of
$s^n(\bblanc)$ we find a zone of computation of level $m$, then in this line we
have $4^{n-m}$ zones of computation of level $m$. Moreover for each computation
box located at the coordinate $(i,j)$, the next computation box in the same
column above the current one is
separated by $2^m-1$ communication boxes, so it is located at the coordinate
$(i,j+2^m)$, and this computation box is in a zone of computation of level $m$
at the same place that the box at the position $(i,j)$. There is the same
phenomena if we look down. The set of zones of computation of the same size
$2^m$ on a vertical line is called a {\em strip of computation} of size $2^m$. 

For any other symbol $a\in\mathcal{G}_1$, the description is the same except for
the bottom row.

\begin{remark}
Note that it is possible to have a $\bleft$ symbol on a row with no $\bright$
symbol on its right -- that is to say an infinite computation zone. In this case
the computation zone has an infinite level.
\end{remark}

\begin{fact}\label{fact.Tgrid}
Consider $x\in\T_{\texttt{Grid}}$ and $C_n$ a computation zone of level $n$ of
$x$. We assume that $C_n$ appears in the i$^\text{th}$ row of $x$, and we denote
this row by $x_{\Z\times \{i\}}$. Then the following properties hold -- remember
that in this section we only consider the $\mathcal{G}_1$ part of the subshift
$\T_{\texttt{Grid}}$:
\begin{enumerate}
\item $\mathcal{C}_n$ contains $2^n$ computation boxes, separated by
communication boxes;
\item on $x_{\Z\times \{i\}}$ there are only computation zones of level $n$,
separated by
$2^{2n-1}$ communication boxes;
\item the row $x_{\Z\times \{i\}}$ is repeated vertically every $2^n$ rows, that
is to say
$x_{\Z\times \{i\}}=x_{\Z\times \{i+k\times2^n\}}$ for any integer $k\in\Z$;
\item vertically, between every pair of consecutive computation boxes of $x_{\Z\times \{i\}}$ and
$x_{\Z\times \{i+k\times2^n\}}$,
there are only $2^n-1$ communication boxes.
\end{enumerate}
\end{fact}

We now explain how it is possible for two computation boxes of the same
computation zone and for two adjacent strips to communicate.

\paragraph{Communication inside a strip}
Two computation zones in the same computation strip of level $n$ communicate
thanks to communication boxes of the $2^n-1$ intermediate rows (vertical
transfer of information), and inside a same computation zone communication
between computation boxes occurs on the $2*4^{n-1}-2^n$ communication boxes
(horizontal transfer of information). 

\begin{figure}[!ht]
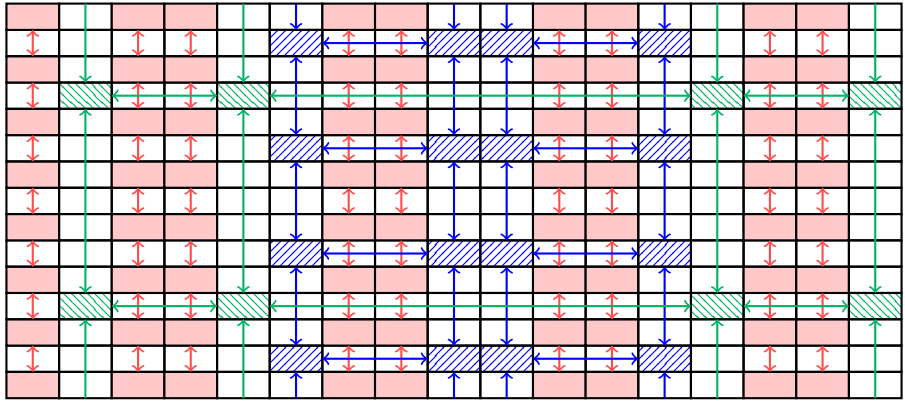
\centering
\include{figure_grid}
\caption{Computation grid with the communication between disconnected parts of
the same computation zone. Computation zones of level 1,2 and 3 are
pictured with three different colors and patterns.
On a given row there are only computation zones of the same level, and there are
$2^n$ rows between two rows with level $n$ computation zones. A row of a strip
of computation of level $n$ is made of $2^n$ boxes arranged into a $2*4^{n-1}$
wide block of boxes. The two ways of communication (horizontal and vertical) are
pictured with arrows whose color corresponds to the level of the computation
zone or strip.}
\label{figure.grid}
\end{figure}

Figure~\ref{figure.grid} represents a computation grid where all zones of
computation of the same size share the same color and are filled with the same
pattern.

\paragraph{Communication between two adjacent strips}
Two strips of computation of same level can also communicate if they are
adjacent -- that is the leftmost computation box of the first strip and the
rightmost computation box of the second strip are only separated by
communication boxes.

\begin{fact}\label{fact.communication}
Communication boxes contain two communication channels -- horizontal and
vertical channels. Thanks to these channels, computation boxes into a same strip
and two adjacent strips can communicate.
\end{fact}

\subsection{Initialization of calculations : the clock}\label{subsection.clock}

The computation strips described in the previous section are restricted in space
but not in time, hence inconsistent configurations of a Turing machine may
appear. To solve this problem, we equip each computation strip with a clock,
that will be reinitialized periodically. At each step of calculation, the clock
is increased and when it is reinitialized, the Turing machine starts a new
calculation.

We use a four elements alphabet $\mathcal{C}=\{ 0,1,\emptyset,\sim\}$ to construct a
sofic subshift $\T_\texttt{Clock}$ obtained by adding finite type rules on
$\product{\T_\texttt{Grid},\mathcal{C}^{\Z^2}}$, where $\T_\texttt{Grid}$ is the
sofic subshift  described in
Section~\ref{subsection.substsofic}. Denote $\pi_{\mathcal{C}}$ the projection on the second coordinate. The clock is actually a finite automaton
that simulates binary addition modulo $2^{2^n}$ on a $2^n$ boxes tape ---
special symbol $\emptyset$ corresponds to the carry in binary addition, and symbol
$\sim$ is used to synchronize adjacent computation zones of same level. To
prevent the appearance of inconsistent states on the clock, we forbid the
patterns $\begin{array}{|c|c|} \hline \emptyset & 0\\ \hline \end{array}$,
$\begin{array}{|c|c|} \hline \emptyset & 1\\ \hline \end{array}$, $\begin{array}{|c|c|}
\hline 0 & \emptyset\\ \hline \end{array}$, $\begin{array}{|c|c|} \hline x & \sim\\
\hline \end{array}$ and $\begin{array}{|c|c|} \hline \sim & x\\ \hline
\end{array}$ where $x\in\{0,1,\emptyset \}$ --- we call this finite type condition
\textbf{Consist}.

\begin{figure}[!ht]
\centering

\begin{tabular}{cc}
\begin{tiny}
\begin{tabular}{cccc}
$\sim$&$\sim$&$\sim$&$\sim$\\
$\emptyset$&$\emptyset$&$\emptyset$&$\emptyset$\\
$1$&$1$&$1$&$1$\\
$1$&$1$&$1$&$\emptyset$\\
$1$&$1$&$0$&$1$\\
$1$&$1$&$\emptyset$&$\emptyset$\\
$1$&$0$&$1$&$1$\\
$1$&$0$&$1$&$\emptyset$\\
$1$&$0$&$0$&$1$\\
$1$&$\emptyset$&$\emptyset$&$\emptyset$\\
$0$&$1$&$1$&$1$\\
$0$&$1$&$1$&$\emptyset$\\
$0$&$1$&$0$&$1$\\
$0$&$1$&$\emptyset$&$\emptyset$\\
$0$&$0$&$1$&$1$\\
$0$&$0$&$1$&$\emptyset$\\
$0$&$0$&$0$&$1$\\
$0$&$0$&$0$&$0$\\
\end{tabular}
\end{tiny}
&{\vbox to 2.5cm{\hbox{\includegraphics[width=13cm]{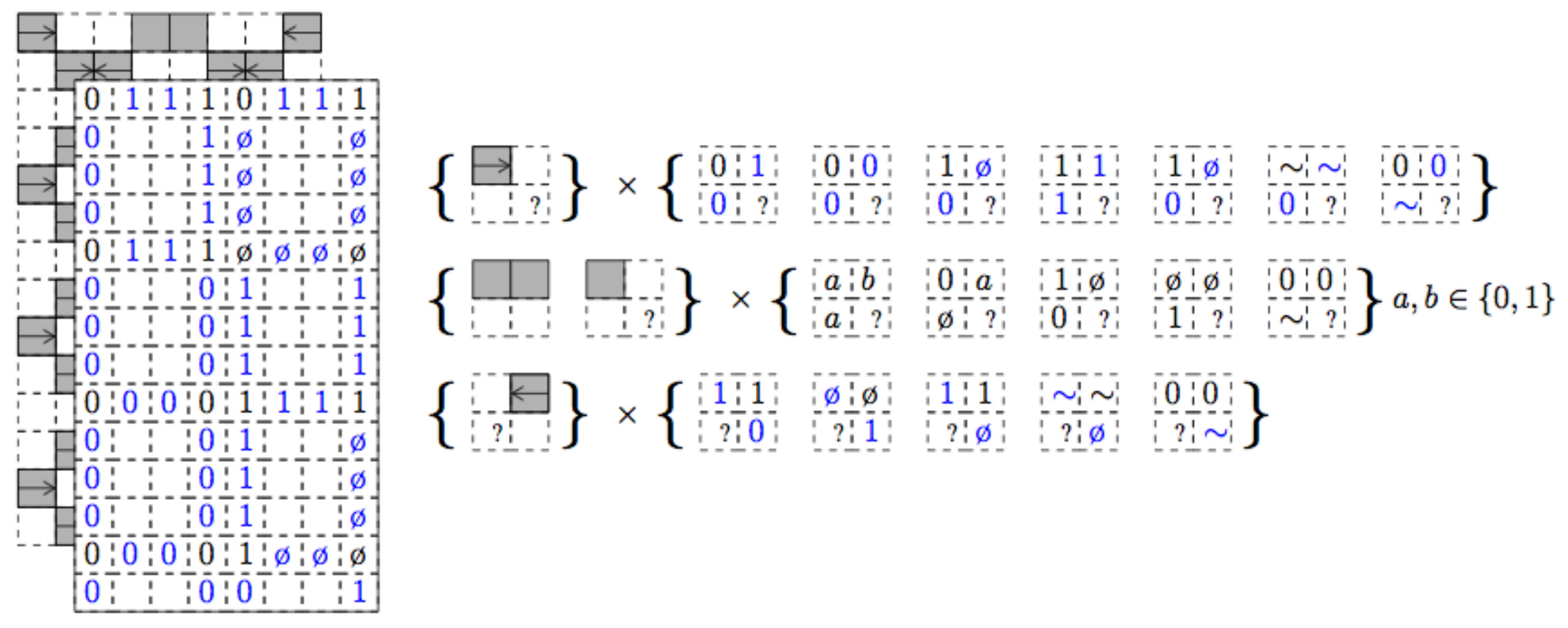}}}}\\
\end{tabular}
\vspace{1cm}
\caption{On the left, an example of the evolution of the clock for a computation
zone of size $2^2$. On the middle the evolution of a part of this clock
on a level 2 computation strip: on the tape are successively written $001\emptyset$,
$0011$, $01\emptyset\emptyset$ and $0101$. And on the right, some of
the finite type conditions \textbf{Count}, represented by the allowed patterns,
added to the sofic subshift $\T_\texttt{Grid}$ to obtain the sofic subshift
$\T_\texttt{Clock}$.}
\label{figure.clock}
\end{figure}

We describe the finite type conditions $\textbf{Count}$ on the alphabet
$\mathcal{G}_1\times \mathcal{C}$ in Figure~\ref{figure.clock}.

The clocks of different computation levels evolve according to the rules
described in Figure~\ref{figure.clock}, and when a symbol $\emptyset$ reaches the left
most computation box $\bleft$, it is reinitialized. Before reinitialization, the
clock passes through the configuration with only $\sim$ symbols on the tape.
Thanks to this configuration, it is possible to synchronize a clock on a strip
of level $n$ with its two neighbours of level $n$. For example the clock for a
computation strip of level $1$ will be $00,01,1\emptyset,11,\emptyset\emptyset,\sim\sim,00,\dots$
Hence a clock for a computation strip of level $n$ is reinitialized after
$2^{2^n}+2$ steps.

To these local rules we add another finite type condition called
\textbf{Synchro}, that ensures that clocks corresponding to computation zones on
the same level are synchronized, that is they are in the same state at every
calculation step -- on a same row, all the clocks are in the same state. This
can be easily done by the following way: a clock is in the configuration
$\sim\dots\sim$ only when its left and right neighbours are in the same
configuration -- a signal carrying symbol $\sim$ is sent through communication
channel between neighbours. We thus obtain a sofic subshift
$$\T_\texttt{Clock}=\ft{\textbf{Count}\cup\textbf{Consist}\cup
\textbf{Synchro}}{\product{\T_\texttt{Grid},\mathcal{C}^{\Z^2}}}$$
in which every computation strip of $\T_\texttt{Grid}$ is now equipped with a
clock. Note that we do not impose clocks for different levels of computation
zones to be somehow synchronized.

\begin{fact}
Consider the sofic subshift $\T_\texttt{Clock}$, in the interior of a strip of
computation of level $n$ which is of size $2^n$, the clock is initialized every
$2^{2^n}+2$ on computation zone of level $n$.
\end{fact}

\subsection{A sofic subshift to describe Turing machines
behaviour}\label{subsection.sftMt}

We are going to use the subshift $\T_{\texttt{Clock}}$ constructed in
Section~\ref{subsection.clock} to construct a sofic subshift where the
computation of $\TM$ in a space $2^n$ holds on each strip of computation of size
$2^n$, for all $n\in\N^{\ast}$. We want to apply the rules of $P_{\TM}$ to
adjacent computation boxes that may be separated by a sequel of communication
boxes. As explained in Section~\ref{subsection.communication} information may be
transferred through communication boxes horizontally and vertically. The space
of computation of $\TM$ is restricted by $\bleft$ on the left and by $\bright $
on the right. We start again with the sofic subshift $\T_\texttt{Clock}$ defined
in Section~\ref{subsection.clock}, into the product subshift
$\product{\T_{\texttt{Clock}},\tilde{\mathcal{A}}^{\Z^2}}$ where
$\tilde{\mathcal{A}}=\mathcal{A}_{\TM}\cup(\mathcal{A}_{\TM}\times
\mathcal{A}_{\TM}\times
\mathcal{A}_{\TM})$. A symbol in $\tilde{\A}$ may be either a symbol of
$\mathcal{A}_{\TM}$ inside a computation box or three symbols of
$\mathcal{A}_{\TM}$ transferred -- horizontally for the first and the second and vertically for
the third -- through a communication box. We have defined $\pi_{\mathcal{G}_1}$, $\pi_{\mathcal{G}_2}$ and $\pi_{\mathcal{C}}$ respectively the projections on $\mathcal{G}_1$, $\mathcal{G}_2$ and $\mathcal{C}$ in the first coordinate of $\product{\T_{\texttt{Clock}},\tilde{\mathcal{A}}^{\Z^2}}$. Moreover denote $\pi_{\tilde{\A}}$ the projection on the second coordinate of $\product{\T_{\texttt{Clock}},\tilde{\mathcal{A}}^{\Z^2}}$, if we are in a communication box, we can write $\pi_{\tilde{\A}_1}$, $\pi_{\tilde{\A}_2}$ and $\pi_{\tilde{\A}_3}$ respecively for the first, second and third coordinate of $\mathcal{A}_{\TM}\times\mathcal{A}_{\TM}\times
\mathcal{A}_{\TM}$. 

To the sofic-subshift
$\product{\T_{\texttt{Clock}},\tilde{\mathcal{A}}^{\Z^2}}$, we add the following
finite conditions, the support of all forbidden patterns have the following form: $$\begin{array}{c|c|c}
\cline{2-2}&a&\\ \hline \multicolumn{1}{|c|}{b}&c&\multicolumn{1}{c|}{d} \\ \hline&e&\\ \cline{2-2}\end{array} \qquad\textrm{ with }a,b,c,d,e\in\mathcal{G}_1\times\mathcal{G}_2\times\mathcal{C}\times\tilde{\A}$$
The conditions are:
\begin{itemize}
\item if the center box corresponds to a communication box in
$\T_{\texttt{Clock}}$, that is to say $\pi_{\mathcal{G}_1}(c)=\bblanc$, one uses conditions \textbf{Transfer}: the first and second coordinates are constant along the central row, and the third coordinate is constant along the central comlumn -- more precisely $\pi_{\tilde{\A}_1}(b)=\pi_{\tilde{\A}_1}(c)=\pi_{\tilde{\A}_1}(d)$, $\pi_{\tilde{\A}_2}(b)=\pi_{\tilde{\A}_2}(c)=\pi_{\tilde{\A}_2}(d)$ and $\pi_{\tilde{\A}_3}(a)=\pi_{\tilde{\A}_3}(c)=\pi_{\tilde{\A}_3}(e)$, these conditions hold if all boxes in the neighborhood are communication boxes, in fact, if there is a computation box, we just use the projection $\pi_{\tilde{\A}}$;
\item if the center box corresponds to a computation box in
$\T_{\texttt{Clock}}$, that is to say $\pi_{\mathcal{G}_1}(c)\in\left\{\bnoir,\bleft ,\bright\right\}$, one uses one of the followings conditions:
\begin{itemize}
\item conditions \textbf{Init}: when the clock is in a initial state, there is
the blank symbol $\sharp$ on each box and the tape is in the initial state on
the left computation box $\bleft$ -- more precisely 
\begin{itemize}
\item if $\pi_{\mathcal{C}}(c)=\sim$ and $\pi_{\mathcal{G}_1}(c)=\bleft$ then $\pi_{\tilde{\A}}(c)=\pi_{\tilde{\A}_1}(d)=\pi_{\tilde{\A}_2}(b)=\pi_{\tilde{\A}_3}(a)=(q_0,\sharp)$,
\item if $\pi_{\mathcal{C}}(c)=\sim$ and $\pi_{\mathcal{G}_1}(c)\in\{\bnoir,\bright\}$ then $\pi_{\tilde{\A}}(c)=\pi_{\tilde{\A}_1}(d)=\pi_{\tilde{\A}_2}(b)=\pi_{\tilde{\A}_3}(a)=\sharp$;
\end{itemize}
\item conditions \textbf{Comp}: we use the rules described in
$\mathcal{P}_{\TM}$ if the clock is not in the initial state -- more precisely \begin{itemize}
\item if $\pi_{\mathcal{C}}(c)\ne\sim$ and $\pi_{\mathcal{G}_1}(c)=\bnoir$ then 
$$ \begin{array}{c|c|c}
\cline{2-2}&\pi_{\tilde{\A}_3}(a)&\\ \hline \multicolumn{1}{|c|}{\pi_{\tilde{\A}_1}(b)}&\pi_{\tilde{\A}}(c)&\multicolumn{1}{c|}{\pi_{\tilde{\A}_2}(d)}\\ \hline\end{array} \in P_{\TM}, \qquad \pi_{\tilde{\A}}(c)=\pi_{\tilde{\A}_2}(b)=\pi_{\tilde{\A}_1}(d) \textrm{ and } \pi_{\tilde{\A}}(c)=\pi_{\tilde{\A}_3}(e),$$
\item if $\pi_{\mathcal{C}}(c)\ne\sim$, $\pi_{\mathcal{G}_1}(c)=\bleft$ and the third coordinate of $\delta(\pi_{\tilde{\A}}(c))$ is different from $\leftarrow$, that is to say the transition function of the Turing machine does not move the head toward the left, then 
$$ \begin{array}{c|c|c}
\cline{2-2}&\pi_{\tilde{\A}_3}(a)&\\ \hline \multicolumn{1}{|c|}{\sharp}&\pi_{\tilde{\A}}(c)&\multicolumn{1}{c|}{\pi_{\tilde{\A}_2}(d)}\\ \hline\end{array} \in P_{\TM}, \qquad \pi_{\tilde{\A}}(c)=\pi_{\tilde{\A}_2}(b)=\pi_{\tilde{\A}_1}(d) \textrm{ and } \pi_{\tilde{\A}}(c)=\pi_{\tilde{\A}_3}(e),$$
if $\pi_{\mathcal{C}}(c)\ne\sim$, $\pi_{\mathcal{G}_1}(c)=\bright$ and the third coordinate of $\delta(\pi_{\tilde{\A}}(c))$ is different from $\rightarrow$, that is to say the transition function of the Turing machine does not move the head toward the right, then 
$$ \begin{array}{c|c|c}
\cline{2-2}&\pi_{\tilde{\A}_3}(a)&\\ \hline \multicolumn{1}{|c|}{\pi_{\tilde{\A}_1}(b)}&\pi_{\tilde{\A}}(c)&\multicolumn{1}{c|}{\sharp}\\ \hline\end{array} \in P_{\TM}, \qquad \pi_{\tilde{\A}}(c)=\pi_{\tilde{\A}_2}(b)=\pi_{\tilde{\A}_1}(d) \textrm{ and } \pi_{\tilde{\A}}(c)=\pi_{\tilde{\A}_3}(e);$$
\end{itemize}
\item conditions \textbf{Bound}: if the head wants to go to the left of the
computation box $\bleft$ or to the right of the computation box $\bright$, the head reaches a special state and the computation
continues in an infinite loop until the computation is initiated by the clock -- more precisely 
\begin{itemize}
\item if $\pi_{\mathcal{C}}(c)\ne\sim$, $\pi_{\mathcal{G}_1}(c)=\bleft$ and the third coordinate of $\delta(\pi_{\tilde{\A}}(c))$ is $\leftarrow$, then 
$$ \begin{array}{|c|c}
\cline{1-1}q_{\texttt{Wait}}&\\ \hline \pi_{\tilde{\A}}(c)&\multicolumn{1}{c|}{\pi_{\tilde{\A}_2}(d)}\\ \hline\end{array}, \qquad \pi_{\tilde{\A}}(c)=\pi_{\tilde{\A}_2}(b)=\pi_{\tilde{\A}_1}(d) \textrm{ and } \pi_{\tilde{\A}}(c)=\pi_{\tilde{\A}_3}(e);$$
\item if $\pi_{\mathcal{C}}(c)\ne\sim$, $\pi_{\mathcal{G}_1}(c)=\bright$ and the third coordinate of $\delta(\pi_{\tilde{\A}}(c))$ is $\rightarrow$, then 
$$ \begin{array}{c|c|}
\cline{2-2}&q_{\texttt{Wait}}\\ \hline \multicolumn{1}{|c|}{\pi_{\tilde{\A}_1}(b)}&\pi_{\tilde{\A}}(c)\\ \hline\end{array}, \qquad \pi_{\tilde{\A}}(c)=\pi_{\tilde{\A}_2}(b)=\pi_{\tilde{\A}_1}(d) \textrm{ and } \pi_{\tilde{\A}}(c)=\pi_{\tilde{\A}_3}(e);$$
\item  if $\pi_{\mathcal{C}}(c)\ne\sim$, $\pi_{\mathcal{G}_1}(c)\in\{\bnoir,\bleft,\bright\}$ and $\pi_{\tilde{\A}}=q_{\texttt{Wait}}$, then 
$$  \pi_{\tilde{\A}}(c)=\pi_{\tilde{\A}_3}(a)=\pi_{\tilde{\A}_3}(e)=\pi_{\tilde{\A}_2}(b)=\pi_{\tilde{\A}_1}(d)=q_{\texttt{Wait}},$$
\item if $\pi_{\mathcal{C}}(c)\ne\sim$ and $\pi_{\tilde{\A}_1}(b)=q_{\texttt{Wait}}$ or $\pi_{\tilde{\A}_2}(d)=q_{\texttt{Wait}}$ then $\pi_{\tilde{\A}}(c)=q_{\texttt{Wait}}$.
\end{itemize}
\end{itemize}
\end{itemize}

Define the sofic subshift $\T_{\TM}$:
$$\T_{\TM} =
\ft{\textbf{Transfer}\cup\textbf{Init}\cup\textbf{Comp}\cup\textbf{Bound}}{
\product{\T_{\texttt{Clock}},\tilde{\mathcal{A}}^{\Z^2}}}.$$

For more convenience, we gather the local rules \textbf{Transfer},
\textbf{Init}, \textbf{Comp} and \textbf{Bound} in $\textbf{Work}_\TM$, and the
construction is summed up by: $\T_{\TM} =
\ft{\textbf{Work}_{\TM}}{\product{\T_{\texttt{Grid}},\mathcal{A}_{Comp(\TM)}^{
\Z^2}}}$ for any Turing machine $\TM$.

On each strip of computation appears parts of the space time diagram of the
calculation of $\TM$ on the empty word. Each part of these space time diagrams are limited in space by the size of the  strip of computation and the number of steps is bounded exponentially by the length of the strip. Thus we can find in $\T_{\TM}$ arbitrary large part of space time diagram of $\TM$.

\begin{fact}\label{fact.TMtime}
The subshift $\T_{\TM}$ contains all calculations of the Turing machine $\TM$ on
space time diagram of size $2^n\times \left(2^{2^n}+2\right)$ -- $2^n$ boxes tape and $2^{2^n}+2$ steps of
calculation -- starting with an empty entry word.
\end{fact}

\begin{example}

In this example the Turing machine $\TM_{\texttt{ex}}$ starts its enumeration
with the word $ab$. The picture describes how a run is coded on a computation
grid. If one only considers computation boxes of level $2$, they form a three by
four computation zone (three steps of calculation on a four boxes tape).

\begin{figure}[!ht]
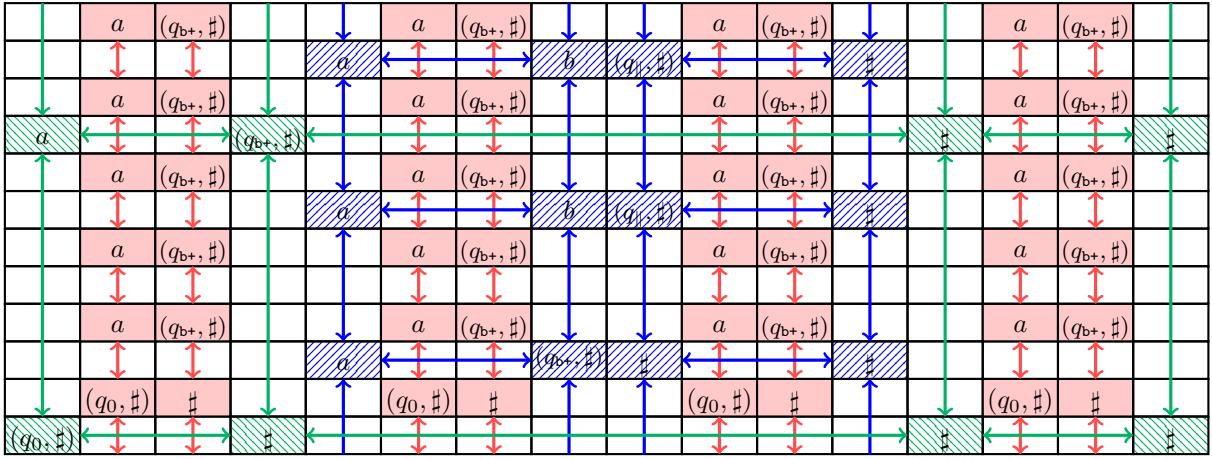
\centering
\include{figure_calcul_fractionne}
\caption{Calculation of a Turing machine on a computation grid with computation
zones of levels 1, 2 and 3. Remark that each $\uparrow$ or $\leftrightarrow$
arrow actually carries a symbol, but for more readability they are not pictured
here. For the same reason the clock is also omitted.}
\label{figure.calculfractionne}
\end{figure}
\end{example}


\subsection{Communication channels between Turing machine of different
levels}\label{subsection.channels}

In the sequel computation strips will need to communicate. For two strips of
the same level communication it is easy since between two zones of computation of
adjacent strips of level $n$, there are only communication boxes. Then one bit
of information can be exchanged between two adjacent strips of level $n$ at each
step of calculation (see Section~\ref{subsection.computationzones}). But if the
two strips are not of the same level the problem is not as simple. We present in
this section a communication grid that allows a strip of level $n$ to
communicate with a strip of level $n-1$ and a strip of level $n+1$. This
communication grid is based on the $\mathcal{G}_2$ part of the subshift
$\T_{\texttt{Grid}}$.

The lines obtained with the alphabet $\mathcal{G}_2$ are called {\em
communication lines}. Communication between computation zones of different
levels are made through these lines. Under the action of $s_2$, communication
lines form rectangles. The two rectangles obtained after $n$ iteration of an element of $\mathcal{G}_2$ are called  
{\em communication rectangles of level $n$}. Each rectangle of level $n$
intersects two rectangles of level $n-1$ and it is intersected by a rectangle of
level $n+1$.

If we consider a border computation box $\bleft$ (resp. $\bright$) in a computation zone of level $n$, it is inside a communication rectangle of level $n$. Thus if we go horizontally on the left (resp. the right) of this box we meet the left border (resp. the right border) of this rectangle. On the bottom and top lines of this
rectangle, we encounter two border computation boxes ($\bleft$ and $\bright$)
which are in two different computation zones of level $n-1$.

\begin{figure}[!ht]\centering
\includegraphics[width=16cm]{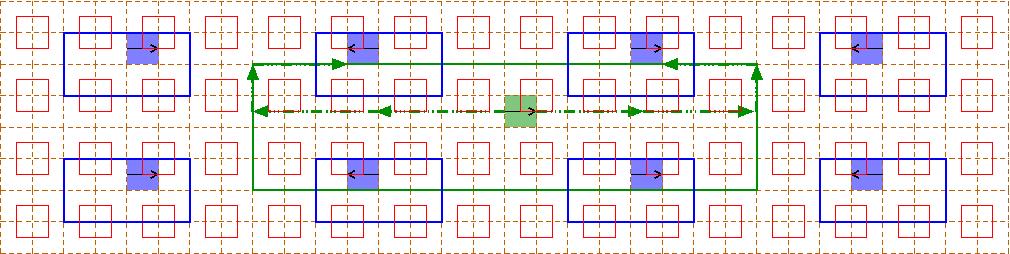}

\caption{A computation grid with communication lines. The computation zone of
level 3 communicates with level 2 computation zones it contains. This
communication is made through the level 3 communication rectangle inside which
the left border computation box is. Symmetrically, one can imagine that the
right border communication box communicates with two other level 2 computation
zones, this in not pictured here.}
\label{figure.communication}
\end{figure}

By local rules it is possible to construct \emph{communications channel of level
$n$}, that start from each border computation box ($\bleft$ or $\bright$) of
level $n$. The channel of communication goes on horizontally on the right and
left branches until it meets the right or left border of a communication
rectangle which is necessary of level $n$. Then the channel goes up and follows
the border of the rectangle until it meets a border computation box. This box is
necessarily of level $n-1$. Thus a computation zone can communicate with the
four computation zones of the previous level which are included in itself (see
Figure~\ref{figure.communication}). These channels are used in
Section~\ref{subsection.mt2} to ensure communication between computation zones
of different levels. We remark that zones of a level $n$ repeat vertically with
half the frequency of level $n-1$ zones. Therefore half the level $n-1$ zones do
not  incoming path from higher zones.

\begin{fact}\label{fact.chanelcommunication}
For any computation strip of level $n$, there are two communication channels
starting from each border computation box $\bleft$ or $\bright$ of level $n$ and
ending at a border computation box of level $n-1$ -- one $\bleft$ and one
$\bright$. Starting from a computation zone of level $n$, the four computation
strips of level $n-1$ associated can be reached by this way.
\end{fact}




\section{Proof of the main theorem}\label{section.MainResult}


The ideas of the proof of the main result of this article were presented in the
Introduction. We give here technical details that rely on constructions
presented in the previous sections. We want to prove the following result.

\begin{theorem}\label{theorem.mainresult}
Any effective subshift of dimension $d$ can be obtained with factor and projective
subaction operations from a subshift of finite type of dimension $d+1$.
\end{theorem}

Thanks to the formalism of Section~\ref{section.subshifts} and since $\reshift$ is stable under $\Factor$ and $\SA$ operations, we rewrite it:
$$\mathcal{C}l_{\Factor,\SA}(\SFT\cap \shift_{d+1})\cap\shift_{\leq
d}=\reshift\cap\shift_{\leq d}.$$

This result improves Hochman's~\cite{hochman2007drp} since our construction
decreases the dimension.

We here prove this statement in the particular case $d=1$, but the proof can be
easily extended to any dimension. Let $\Sigma$ be a one dimensional effective
subshift, defined on an alphabet $\A_{\Sigma}$. 

\subsection{Construction of the four layers of SFT}\label{subsection.sketch}

We start with the two-dimensional fullshift ${\A_\Sigma}^{\Z^2}$ with a spatial
extension operation, and thanks to factor, product and finite type operations we
construct a sofic subshift $\T_{\texttt{Final}}$ such that after factor and
projective subaction we obtain $\Sigma$. To do that, we eliminate configurations
$x$ such that $x_{\Z\times\{0\}}$ contains a forbidden word of $\Sigma$. Then
the projective subaction that consists in only keeping the first coordinate of a
two-dimensional configuration $x$ gives the subshift $\Sigma$.

To resume the two-dimensional sofic subshift is made of four layers that are
glued together thanks to product operations:
\begin{itemize}
 \item first layer contains $\A_{\Sigma}^{\Z^2}$ and all horizontal lines are
identical by finite condition $\textbf{Align}$, the other layers force the
horizontal line to be an element of the effective subshift $\Sigma$, thus this
subshift can be obtained after projective subaction (to keep horizontal line)
and factor (to keep the first layer);
 \item layer 2 contains the computation zones for Turing machines equipped with
the clock (this construction is described in
Section~\ref{section.computationzones}), that will be used by both machines
$\TM_{\texttt{Forbid}}$ and $\TM_\texttt{Search}$; but also the communication
channels that will be used by the same machines to send requests (see
Sections~\ref{subsection.mt1} and~\ref{subsection.mt2});
 \item layer 3 is devoted to Turing machines $\TM_{\texttt{Forbid}}$, and
communication with the Turing machines $\TM_\texttt{Search}$ (this part is
described in Section~\ref{subsection.mt1});
 \item layer 4 is devoted to Turing machines $\TM_\texttt{Search}$ and internal
communication between these machines (see Section~\ref{subsection.mt2}).
\end{itemize}
Of course each of this layer depends on the others (for example layer 3 uses
computation zones given by layer 2), and the dependences are coded thanks to
finite type operations.

\begin{figure}[!ht]\begin{center}
\includegraphics[width=7.5cm]{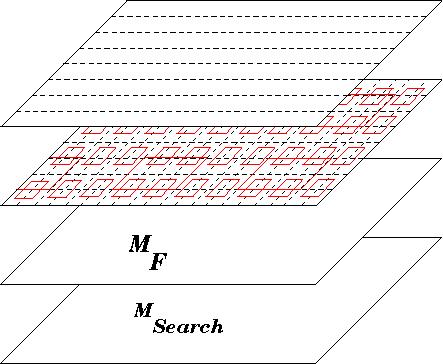}
\end{center}
\caption{Four layers in the final construction.}
\label{figure.level}
\end{figure}

\subsection{Addresses in a strip}\label{subsection.addresses}
 
Since, on the first layer, each column is formed by one letter of $\A_{\Sigma}$,
to check a word in an horizontal configuration, it is sufficient to check the
first layer in the corresponding columns.

Let $C_n$ be a computation zone of level $n$ of an element $x\in\T_{\texttt{Grid}}$ and let $S_n$ be the computation strip associated. By Fact~\ref{fact.UniqueDerivation}, there exists an unique $i\in[0,4^n-1]\times[0,2^n-1]$ and an unique $y\in\T_{\texttt{Grid}}$ such that $s_{\texttt{Grid}}^n(y)=\s^i(x)$ so there exists an unique $(j_1,j_2)\in\Z^2$ such that $C_n\sqsubset s_{\texttt{Grid}}^n(y_{(j_1,j_2)})$. One has $S_n\sqsubset \sigma^{-i}(s_{\texttt{Grid}}^n(y_{\{j_1\}\times\Z\}})) \sqsubset x$, the strip $\sigma^{-i}(s_{\texttt{Grid}}^n(y_{\{j_1\}\times\Z\}}))$ is the \emph{dependency strip} associated with the computation strip $S_n$ 

In $\T_{\texttt{Grid}}$ the tape of a Turing machine in a strip of level $n$ is fractured. Thus a Turing machine of level $n$
cannot view all columns which are in its associated dependency strip. To get
this information, this Turing machine communicates with a Turing machine of lower
level (see Section~\ref{subsection.mt2}) but both machines need to precisely
identify a column.  

\begin{figure}[!ht]\begin{center}
\begin{tikzpicture}[scale=0.9]
 \draw (0,0) rectangle (16,2);
 \foreach \i in {0,...,3}{
 \draw (4*\i,0) rectangle (4*\i+4,-1);
 }
  \foreach \i in {0,...,15}{
 \draw (\i,-1) rectangle (\i+1,-1.5);
 }
   \foreach \i in {0,...,63}{
 \draw (0.25*\i,-1.5) rectangle (0.25*\i+0.25,-1.75);
 }
  \foreach \i in {0,...,3}{
 \draw (2+4*\i,0) node[above]{$\i$};
  \foreach \j in {0,...,3}{
  \draw (4*\i+\j+0.5,-1) node[above]{$\j$}; 
 } }
 \draw[fill=gray!70] (2,-1.5) rectangle (2.25,-1.75);
  \draw[fill] (11.25,-1.5) rectangle (11.5,-1.75);

\draw (16,1) node[right] {$s^n_1$};
\draw (16,-0.5) node[right] {$s^n_2$};
\draw (16,-1.25) node[right] {$s^n_3$};
\end{tikzpicture}
\end{center}
\caption{Addresses of two boxes inside a dependency strip associated with a
computation zone of level~3. The address of the column of the black box is $231$
and for the grey box, the address of the column is $020$.}
\label{figure.address}
\end{figure}
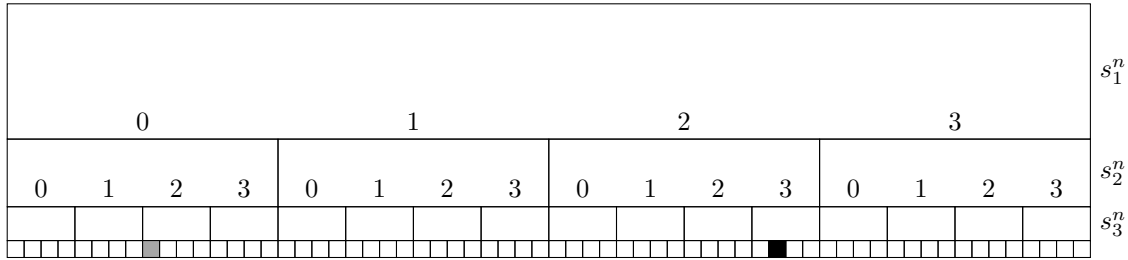

Given a dependency strip associated with a computation strip of level $n$, it is
possible to describe the coordinate relative to this strip of any
column of the dependency strip by an address which contains $n$ letters
in a four elements alphabet. Each $s_1^n(a)$ is horizontally decomposed into
four (possibly different) $s_1^{n-1}(b)$ where $a,b\in\mathcal{G}_1$. The first
letter of the address indicates in which of these dependency stripes of size
$n-1$ the column is located. By iteration of this process the position of a
column is
exactly given with $n$ letters (see Figure~\ref{figure.address}).

\begin{fact}\label{fact.address}
For every dependency strip associated with a computation strip of level $n$, it
is possible to describe the position of any column by an address of length $n$
on a
four elements alphabet.
\end{fact}


\subsection{Responsibility zones}\label{subsection.responsibility}

On each computation zone a Turing machine makes calculations. The Turing machine
$\TM_{\texttt{Forbid}}$ described more precisely in Section~\ref{subsection.mt1}
enumerates patterns and then checks that these patterns never appear. Since it
takes an infinite number of steps of calculation to check that one pattern does
not appear in the entire configuration, each Turing machine
$\TM_{\texttt{Forbid}}$ only checks a finite zone. The finite zone in which the
machine ensures that no forbidden pattern it produces appears is called the
\emph{responsibility zone} of the machine.

\begin{figure}[!ht]\centering
\includegraphics[width=18cm]{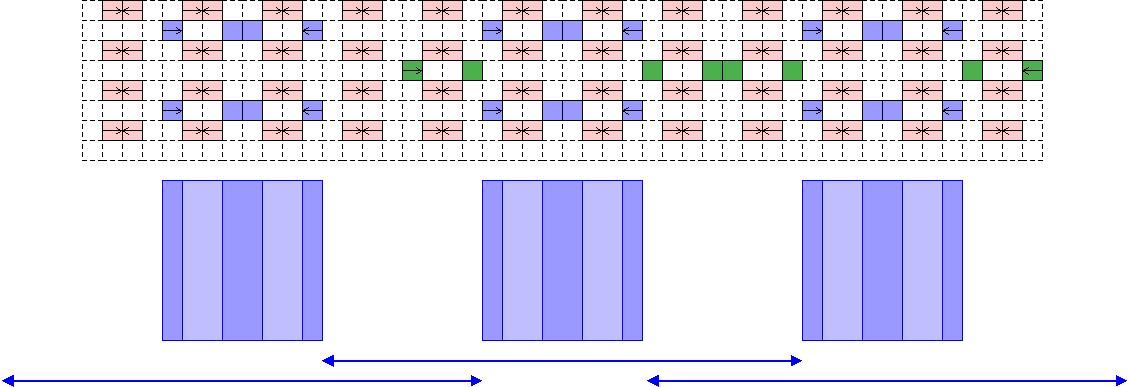}
\caption{Responsibility zones for strips of level 2. These zones are 24 boxes
wide and overlap on 8 boxes. The responsibility zone of the center strip starts
at the end of the strip on its left and ends at the beginning of the strip on
its right.}
\label{figure.responsibilityzone}
\end{figure}

We thus associate a responsibility zone with each strip of computation. For a
strip of level $n$ this responsibility zone is $3*(2*4^{n-1})=6*4^{n-1}$ wide
and centered on the strip (see Figure~\ref{figure.responsibilityzone}), so that
the responsibility zone of a strip starts at the end of the strip of same level
on its left and ends at the beginning of the strip of same level on its right.

Responsibility zones defined in this way overlap: two adjacent responsibility
zones of same level $n$ share $2*4^{n-1}$ boxes. These overlappings are
essential: if they did not exist, one can imagine that a forbidden pattern not
entirely included in any responsibility zone would not be detected. Moreover the
non bounded size of overlappings ensures that any pattern is inside an infinite
number of responsibility zones of increasing levels.


\subsection{Generation and detection of forbidden patterns by
$\TM_{\texttt{Forbid}}$}\label{subsection.mt1}

Since $\Sigma$ is recursively enumerable, there exists a Turing machine that
enumerates the forbidden patterns of $\Sigma$. We here describe a modified
version of this Turing machine that also checks that no forbidden pattern
appears inside its responsibility zone, on the first level of the construction
${\A_\Sigma}^{\Z^2}$. Computation zones are not connected (see
Figure~\ref{figure.calculfractionne}), so a calculation of
$\TM_{\texttt{Forbid}}$ on a strip of computation of size $2^n$ cannot access
entirely its responsibility zone. The machine $\TM_{\texttt{Forbid}}$ needs the help of a
second Turing machine $\TM_\texttt{Search}$ to obtain the patterns of
${\A_\Sigma}^\Z$ written in its responsibility zone. The behaviour of
$\TM_{\texttt{Forbid}}$ is the following: it enumerates as many forbidden
patterns as the size of the computation zone allows, and each time such a pattern
is generated, $\TM_{\texttt{Forbid}}$ checks that it does not appear in its
responsibility zone. 

\paragraph{Tapes of $\TM_{\texttt{Forbid}}$} The machine $\TM_{\texttt{Forbid}}$ uses three tapes:
\begin{itemize}
 \item the first tape is the \emph{calculation tape};
 \item the second tape is a \emph{writing tape}, where the forbidden patterns
are successively written;
 \item the last tape is the \emph{communication tape} and contains successively
the addresses of letters from alphabet $\A_\Sigma$ needed by
$\TM_{\texttt{Forbid}}$ to check no forbidden pattern appears inside its
responsibility zone; $\TM_{\texttt{Forbid}}$ waits for the required
$\TM_\texttt{Search}$ machine of its neighbourhood (left, middle or right
machine) to be available, then sends it the address of the letter it wants to
access (see Section~\ref{subsection.mt2}).
\end{itemize}

\paragraph{Detection of the size of the responsibility zone associated} First, the Turing machine $\TM_{\texttt{Forbid}}$ detects the size of the
computation zone between $\bleft$ and $\bright$. Thus, $\TM_\texttt{Forbid}$
knows the size of its responsibility zone. This can be in linear time according
to the size of the computation zone considered. 

\paragraph{Enumeration of forbidden patterns} Then, $\TM_{\texttt{Forbid}}$ enumerates forbidden patterns and each time it encounter one, it checks if this forbidden pattern appears in the associated responsibility zone before to enumerate the following one.  

\paragraph{Check of the responsibility zone} Assume that the machine $\TM_{\texttt{Forbid}}$ has written on its
writing tape a forbidden pattern $f=f_0 f_1\dots f_{k-1}$. Assume that $\TM_{\texttt{Forbid}}$ must check a responsibility zone of level $n$ denoted $a_0a_1\dots a_{6*4^{n-1}-1}$. It asks
$\TM_\texttt{Search}$ for the first letter in its responsibility zone $a_0$
(the principe of a request is explained in Section~\ref{subsection.mt2}), and compares
it with $f_0$. If the letters coincide, then it is still possible that $f$
appears in position $0$ in the responsibility zone, so the comparison of the two
patterns $f$ and $a_0\dots a_k$ continues. If $f_0\neq a_0$ then we are sure
that $f$ does not appear at this location. If $f=a_0\dots a_k$, the Turing machine $\M_{\texttt{Forbid}}$ stops its computation and enter in a state which says that a forbidden patterns appears in the checked configuration. This state will be forbidden in the final subshift of finite type. When the word $a_0\dots a_{k-1}$ is checked, $\TM_{\texttt{Forbid}}$ continues the comparison with $a_1\dots
a_{k}$, $\dots$, $a_{6*4^{n-1}-k-1}\dots a_{6*4^{n-1}-1}$. At most, to check if $f$ appears in the responsibility zone of level $n$ , $\TM_{\texttt{Forbid}}$ takes $6*4^{n-1}*k*t(n)$  where $t(n)$ is the time takes by 
$\TM_\texttt{Search}$ to answer a request of $\TM_{\texttt{Forbid}}$; the time $t(n)$ is estimated in Section~\ref{subsection.mt2}. 

\begin{center}
\begin{minipage}{9cm}
\begin{figure}[H]
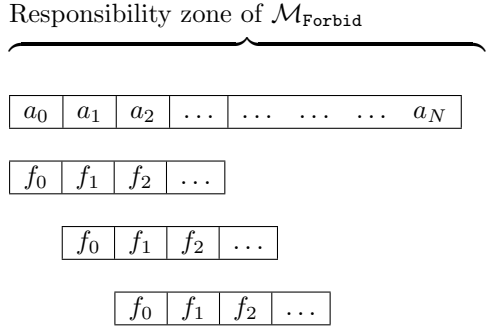
\label{figure.patternchecking}
$$
\begin{array}{l}
\text{Responsibility zone of }\TM_{\texttt{Forbid}}\\
\overbrace{~~~~~~~~~~~~~~~~~~~~~~~~~~~~~~~~~~~~~~~~~~~~~~~~~~~~~~}\\
\\
\begin{array}{|c|c|c|c|cccc|c|}
\hline
a_0 & a_1 & a_2 & \dots & \dots & \dots & \dots & a_N \\
\hline
\end{array}

\\
\\
\begin{array}{|c|c|c|c|}
\hline
f_0 & f_1 & f_2 & \dots \\
\hline
\end{array}
\\
\\
\hspace{0.7cm}\begin{array}{|c|c|c|c|}
\hline
f_0 & f_1 & f_2 & \dots \\
\hline
\end{array}
\\
\\
\hspace{1.4cm}\begin{array}{|c|c|c|c|}
\hline
f_0 & f_1 & f_2 & \dots \\
\hline
\end{array}
\end{array}
$$
\caption{When a forbidden pattern $f=f_0 f_1 \dots f_k$ is generated by
$\TM_{\texttt{Forbid}}$, comparisons with the patterns appearing in the
responsibility zone of $\TM_{\texttt{Forbid}}$ are made in parallel.}
\end{figure}
\end{minipage}
\end{center}


\subsection{Scan of the entire responsibility zone by
$\TM_\texttt{Search}$}\label{subsection.mt2}

\begin{figure}[!ht]\centering
\includegraphics[width=15cm]{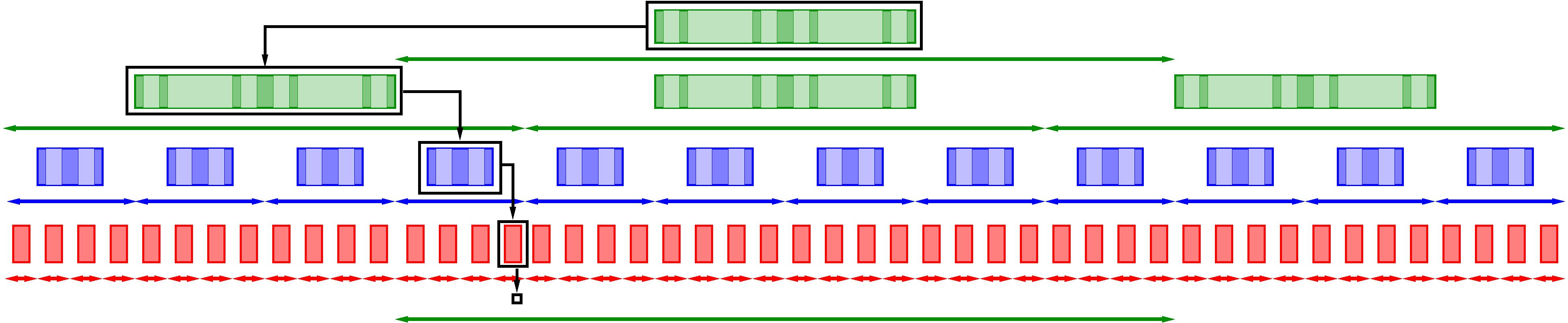}
\caption{An example of request by a $\TM_{\texttt{Forbid}}$ machine of level 3
-- the computation zone on the top of the picture. Depending on the address of
the letter requested, $\TM_{\texttt{Forbid}}$ sends its request to either the
left, center or right $\TM_\texttt{Search}$ machine. On this example the
$\TM_{\texttt{Forbid}}$ machine sends its request to the left
$\TM_\texttt{Search}$ machine of level 3, which transmits it to a
$\TM_\texttt{Search}$ machine of level 2 and finally to a $\TM_\texttt{Search}$
machine of level 1. This last machine can answer the request.}
\label{figure.arbre_recherche}
\end{figure}

The Turing machine $\TM_\texttt{Search}$ is sent a \emph{request} -- that is to say a sequence of symbols which codes the address of a letter inside a responsibility zone of a $\TM_\texttt{Forbid}$ machine -- by
$\TM_{\texttt{Forbid}}$ each time an address is totally written on the
communication tape (the third tape of $\TM_{\texttt{Forbid}}$). The Turing
machine $\TM_\texttt{Search}$ must respond the letter corresponding to the
address inside the responsibility zone, on the first level of the construction
${\A_\Sigma}^{\Z^2}$. Note that the responsibility zone of a $\TM_{\texttt{Forbid}}$ machine of level
$n$ does not exactly match with the communication network of
$\TM_\texttt{Search}$ machines of same level. Actually a $\TM_{\texttt{Forbid}}$
machine shares its responsibility zone with three $\TM_\texttt{Search}$
machines, and depending on the address of the bit requested, the
$\TM_{\texttt{Forbid}}$ sends its request to the appropriate
$\TM_\texttt{Search}$ machine (see Figure~\ref{figure.arbre_recherche} for an
example). 

\paragraph{Tapes of $\TM_{\texttt{Search}}$}The machine $\TM_\texttt{Search}$ of level $n$ uses three
tapes:
\begin{itemize}
\item the first tape is the \emph{calculation tape};
\item the second tape is the \emph{hierarchical request tape}; this is where the bits of an
address transferred by the $\TM_\texttt{Search}$ of level $n+1$ are written.
\item the three last tapes are the \emph{left request tape}, the \emph{center request tape} and the  \emph{right request tape} which correspond to the addresses of the bits asked by the Turing machine $\TM_{\texttt{Forbid}}$ of level $n$ localized respectively to the left, inside and to the right of the communication strip of the machine $\TM_{\texttt{Search}}$ considered.
\end{itemize}

\paragraph{Request sent by $\TM_\texttt{Forbid}$} Each time that an address is written on the communication tape of a Turing machine $\TM_{\texttt{Forbid}}$, this machine sends this request to the corresponding $\TM_\texttt{Search}$ of the same level localized in the same communication strip or in communication strips directly to the left or to the right. $\TM_{\texttt{Forbid}}$  sends one bit composing the address every step
of calculation, so that a level $n$ Turing machine sends a bit every
$2^n$ rows -- if we implement Turing machines in the subshift of finite type described in Section~\ref{subsection.sftMt}. Adjacent strips of same level can communicate by communication channels described in Section~\ref{subsection.computationzones} using the fact that in one row there is only computation zones of same level (see Fact~\ref{fact.Tgrid}). The bits of the address are sent one by one, hence the transfer takes
$2^n*n$ rows since the size of the address of the request is $n$. The request is written on the corresponding request tape. $\TM_{\texttt{Forbid}}$ waits for the answer of the corresponding $\TM_\texttt{Search}$ before to continue the computation. 

\paragraph{Request sent by $\TM_\texttt{Search}$} A Turing machine $\TM_\texttt{Search}$ of level $n\geq 2$ can make a request at one of the four Turing machines $\TM_\texttt{Search}$ of level $n-1$ localized in its dependency. The asking machine sends one bit composing the address every step of calculation, so that a level $n$ Turing machine sends a bit every
$2^n$ rows and thus it takes $2^n*n$ rows to transfer the address of size $n$. The machine $\TM_\texttt{Search}$ of level $n$ uses communication channels described in Fact~\ref{fact.chanelcommunication} to communicate: each border computation box $\bleft$ and $\bright$ is surrounded by a rectangle of the same level $n$ which communicates with border
computation box of the previous level $n-1$. 

\paragraph{Treatment of a request}
A machine $\TM_{\texttt{Search}}$ of level $n$ successively responds to the different request tapes. The address of the request tape considered is copied on the computation tape, and the machine keeps in memory to which request tape it is responding. If the machine
$\TM_\texttt{Search}$ is of level $1$, it directly reads the letter of
$\A_\Sigma$. Otherwise the machine $\TM_\texttt{Search}$ of level $n$
transmits the address to the corresponding machine $\TM_\texttt{Search}$ of level $n-1$: the
first letter of the address indicates which channel $\TM_{\texttt{Search}}$ must be
used to send the continuation of the address, converted into a $n-1$ bits address by erasing the first bit of the address. Then the machine $\TM_\texttt{Search}$ of level $n$ waits for the answer, which is obtained when a machine of level $1$ is reached (see Figure~\ref{figure.arbre_recherche}). This letter must be transferred back until it finds the machine which initially made the request. 

\paragraph{Transfer back of the information} When a Turing machine $\TM_\texttt{Search}$ obtains the bit corresponding to the request, it transfers it by the communication channel to the Turing machine which made the request via the request tapes. This operation is instantaneous for two reasons. First there is just one box of information to transmit. Secondly there is just one information on the channel since the corresponding Turing machine waits for an 
answer. A Turing machine $\TM_\texttt{Search}$ eventually answers the request
of the Turing machine $\TM_\texttt{Forbid}$ of the same computation strip, since
every $\TM_{\texttt{Search}}$ alternately works for
$\TM_{\texttt{Forbid}}$ of same level and higher levels $\TM_\texttt{Search}$ machines. 

\paragraph{Initialization of the computations} When the computation is initialized, it is important not to erase the addresses on the request tapes, because Turing machines of higher levels may be waiting for an answer. Requests are only made toward lower level, so they are answered even if the address does not correspond to a real request. 

Another problem of initialization occurs when a Turing machine makes a request, but is initialized before to obtain its answer. Actually in this case we impose that once the Turing machine is initialized, it waits for the answer to its request from the previous computation before to begin a new one.

\paragraph{Time taken by $\TM_\texttt{Search}$ to answer at a request}

Denote $t(n)$ the time that a machine $\TM_\texttt{Search}$ of level $n$ needs
to answer a request from $\TM_{\texttt{Forbid}}$. Since a machine of level $n$ makes a calculation step
every $2^n$ rows,  a machine $\TM_\texttt{Search}$ of level $n$ needs $2^n*t(n)$ rows
to answer a request from $\TM_{\texttt{Forbid}}$. 

A machine $\TM_\texttt{Search}$ of level $n\geq 2$ needs the help of
a machine $\TM_\texttt{Search}$ of level $n-1$: it transfers one by one the $n-1$
bits of the address, one bit is transferred every $2^n$ rows, this takes $n*2^n$ rows. Then
it waits for the $\TM_\texttt{Search}$ of level $n-1$ answer. It is possible that
this $\TM_\texttt{Search}$ of level $n-1$ is already busy, and the level $n$
machine has to wait -- in the worst case three $\TM_\texttt{Forbid}$ machines of
level $n-1$ are already waiting for an answer. Hence the $\TM_\texttt{Search}$ of
level $n-1$ possibly works for the three neighbouring $\TM_\texttt{Forbid}$ machines of level $n-1$, this takes $3\times
t(n-1)$ steps of calculation, before to work for the $\TM_\texttt{Search}$ of
level $n$, this takes $t(n-1)$ steps of calculation. Thus the number of rows used to answer at a request is given by
$$2^n t(n) \leq n* 2^n + 4* 2^{n-1}*t(n-1).$$
We deduce from the previous inequality that $t(n) \leq n^2 2^n$.

\begin{fact}\label{fact.calculationtime}
All requests of $\TM_\texttt{Forbid}$ of level $n$ are handled by the $\TM_\texttt{Search}$
machine of same level in at most $n^22^n$ steps of calculation for large enough $n$.
\end{fact}

\paragraph{Time taken by $\TM_\texttt{Forbid}$ to check if a forbidden word appear} Assume that a Turing machine $\TM_\texttt{Forbid}$ must check if a word $f$ of size $k$ appears in the responsibility zone associated. According to Section~\ref{subsection.mt1}, this takes $6*4^{n-1}*k*t(n) \leq k*n^2*2^{3n+1}$ steps of calculation.

Let $(f_i)_{i\in\N}$ be the enumeration of forbidden patterns by $\TM_\texttt{Forbid}$. Denote $t(f_0,\dots,f_k)$ the time taken by $\TM_\texttt{Forbid}$ to scan if the words $(f_i)_{i\in[0,k]}$ appear in the responsibility zone associated and denote $t'(f_0,\dots,f_k)$ the time taken by $\TM_\texttt{Forbid}$ to compute the words $(f_i)_{i\in[0,k]}$ without scaning the responsibility zone. Thus, the time taken by a Turing machine $\TM_\texttt{Forbid}$ of level $n$ to scan the words $(f_i)_{i\in[0,k]}$ is given by
$$t(f_0,\dots,f_k)\leq t'(f_0,\dots,f_k)+(k+1)*\max\{|f_i|:i\in[0,k]\}*n^2*2^{3n+1}.$$
Since  $t'(f_0,\dots,f_k)$ does not depend of the level of $\TM_\texttt{Forbid}$ and since by Fact~\ref{fact.TMtime} a machine $\TM_\texttt{Forbid}$ of level $n$ could make $2^{2^n}+2$ steps of calculation, there exists a level $n$ such that all Turing machines of level $n$ check that the words $(f_i)_{i\in[0,k]}$ does not appear in their responsibility zones. 

\begin{fact}\label{fact.timecheck}
For all forbidden word of $\Sigma$, there exists $n\in\N$ such that every turing machine $\TM_\texttt{Forbid}$ of level $n$ checks that the word does not appear in its responsibility zone. 
\end{fact}


\subsection{The final construction}\label{subsection.finalconstruction}

We sum up the construction of the final subshift:
\begin{enumerate}
\item First, we construct the four layers:
\[\T_{\texttt{Level}}=\product{\A^{\Z^2},\T_{\texttt{Grid}},\A^{
\Z^2}_{\textrm{Comp}(\TM_{\texttt{Forbid}})},
\A^{\Z^2}_{\textrm{Comp}(\TM_{\texttt{Search}})}};\]

\item then, we align all the letter of the first layer to obtain the same
configuration horizontally
$\T_{\texttt{Align}}=\ft{\textbf{Align}}{\T_{\texttt{Level}}}$;

\item finally, we include the working of $\TM_{\texttt{Forbid}}$ and
$\TM_{\texttt{Search}}$ thanks to
$\textbf{Work}_{\TM_{\texttt{Forbid}}}\cup\textbf{Work}_{\TM_{\texttt{Search}}}$
and we include the communication between the different layers thanks to 
$\textbf{Com}$. Moreover, we include the condition \textbf{Forbid} which exclude
the configuration when $\TM_{\texttt{Forbid}}$ encounters a forbidden pattern. We
obtain:
$$\T_{\texttt{Final}}=\ft{\textbf{Work}_{\TM_{\texttt{Forbid}}}\cup\textbf{Work}
_{\TM_{\texttt{Search}}}\cup 
\textbf{Com}\cup\textbf{Forbid}}{\T_{\texttt{Align}}}.$$
\end{enumerate}
 
The alphabet of $\T_{\texttt{Final}}$ depends of the Turing machine which enumerates the forbidden patterns of $\Sigma$, it is $O((q.a)^3)$ where $q$ is the number of states and $a$ the cardinal of the alphabet of this Turing machine. Moreover the support of the forbidden patterns of $\T_{\texttt{Final}}$ have the following shape {\vbox to 12pt{\hbox{
\begin{tikzpicture}[scale=0.25]
\draw[dashed] (0,0) rectangle (1,1);
\draw[dashed] (2,0) rectangle (1,1);
\draw[dashed] (2,0) rectangle (3,1);
\draw[dashed] (1,1) rectangle (2,2);
\draw[dashed] (1,0) rectangle (2,-1);
\end{tikzpicture}
}}}.
 \vspace{0.3cm}
 
We denote by $\T$ the subshift $\factor{\pi}{\sa{\Z e_1}{\T_{\texttt{Final}}}}$
where $\pi$ is a morphism that only keeps letters from alphabet $\A_{\Sigma}$
from the first layer. We want to compare $\Sigma$ and $\T$.

\paragraph{Any configuration in $\Sigma$ can be obtained ($\Sigma \subseteq
\T$):}

Let $x\in \Sigma$, by construction of $\T_\text{Final}$ it is easy to construct a
two-dimensional configuration $y$ such that $y\in \T_{\texttt{Final}}$ and
$\pi(y_{|\Z e_1})=x$. 

\paragraph{Any configuration constructed is in $\Sigma$  ($\T \subseteq
\Sigma$):}

Let $x\in\T$, we prove that $x\in \Sigma$. By definition there exists
$y\in\T_\text{Final}$ such that $\pi(y|_{\Z e_1})=x$. It is sufficient to prove
that every word in $x$ is in $\Lang(\Sigma)$. Let $w$ be a word that appears in
$x$. Suppose that $w$ is not in $\Lang(\Sigma)$, by Fact~\ref{fact.timecheck}, there exists
$n\in\N$ such that in any computation strip of level $n$, the word $w$ is
checked in the associated dependency strip. In particular the word $w$ will be
compared with any word of length $|w|$ that appears in $x$. Since $w$ appears in
$x$, there would be a computation strip of level $n$ in which the calculation of
$\TM_{\texttt{Forbid}}$ violates the finite type condition $\textbf{Forbid}$.
This proves the inclusion $\T\subseteq\Sigma$.

\subsection{Effective subshift as sub-action of a two-dimensional sofic}

In fact the previous construction gives a more general result. If we consider
$$\begin{array}{lrll}
\mathbf{SA}_{\Z e_1}:&\pi(\T_{\texttt{Final}})&\longrightarrow&\Sigma\\
&x&\longmapsto &x_{\Z\times\{0\}}
\end{array}$$
it is a continuous bijective map. Indeed, for all $x\in\pi(\T_{\texttt{Final}})$, one has $x_{(i,k)}=x_{(j,k)}$ for all $i,j,k\in\Z$ since by condition $\textbf{Align}$ all columns contain the same symbol. Moreover, $\mathbf{SA}_{\Z e_1}\circ\s^{e_1}=\s_{\Sigma}\circ\mathbf{SA}_{\Z e_1}$, thus $\mathbf{SA}_{\Z e_1}$ realizes a conjugation between the dynamical system $(\pi(\T_{\texttt{Final}}),\s^{e_1})$ and $(\Sigma,\s_{\Sigma})$. We deduce the following theorem:

\begin{theorem}
Any effective subshift of dimension $d$ is conjugate to a sub-action of a sofic subshift of dimension $d+1$.
\end{theorem}

\subsection*{Acknowledgements}
The authors are grateful to Michael Schraudner for useful discussions and important remarks about the redaction. We also
want to thank the anonymous referee for his rigorous and detailed 
review which helped us to clarify the paper and Mike Boyle for some comments about sub-action concepts. 
Moreover, this research is partially supported by projects ANR EMC and ANR SubTile.


\begin{thebibliography}{DRS09}

\bibitem[AS09]{AubrunSablik09}
Nathalie Aubrun and Mathieu Sablik.
\newblock An order on sets of tilings corresponding to an order on languages.
\newblock In {\em 26th International Symposium on Theoretical Aspects of
  Computer Science (STACS 2009)}, volume~3, pages 99--110, 2009.

\bibitem[Bea93]{beal1993cs}
M.P. Beal.
\newblock {\em {Codage Symbolique}}.
\newblock Masson, 1993.

\bibitem[Ber66]{berger1966udp}
R.~Berger.
\newblock {\em {The Undecidability of the Domino Problem}}.
\newblock American Mathematical Society, 1966.

\bibitem[Boy08]{boyle2008}
M.~Boyle.
\newblock {Open problems in symbolic dynamics}.
\newblock {\em Contemporary Mathematics}, 468:69--118, 2008.

\bibitem[Dal74]{Myers1974}
Myers Dale.
\newblock Nonrecursive tilings of the plane. ii.
\newblock {\em The Journal of Symbolic Logic}, 39(2):286--294, 1974.

\bibitem[DLS01]{Durand}
Bruno Durand, Leonid Levin, and Alexander Shen.
\newblock Complex tilings.
\newblock In {\em STOC '01: Proceedings of the thirty-third annual ACM
  symposium on Theory of computing}, pages 732--739, New York, NY, USA, 2001.
  ACM.

\bibitem[DRS08]{Durand2}
Bruno Durand, Andrei~E. Romashchenko, and Alexander Shen.
\newblock Fixed point and aperiodic tilings.
\newblock In {\em Developments in Language Theory}, pages 276--288, 2008.

\bibitem[DRS10]{durand:fixed-point}
Bruno Durand, Andrei~E. Romashchenko, and Alexander Shen.
\newblock Fixed-point tile sets and their applications.
\newblock {\em CoRR abs/0910.2415, \texttt{http://arxiv.org/abs/0910.2415}},
  2010.

\bibitem[Han74]{Hanf}
William Hanf.
\newblock Nonrecursive tilings of the plane. i.
\newblock {\em The Journal of Symbolic Logic}, 39(2):283--285, 1974.

\bibitem[Hed69]{hedlund1969eaa}
GA~Hedlund.
\newblock {Endomorphisms and automorphisms of the shift dynamical system}.
\newblock {\em Theory of Computing Systems}, 3(4):320--375, 1969.

\bibitem[Hoc09]{hochman2007drp}
M.~Hochman.
\newblock {On the Dynamics and Recursive Properties of Multidimensional
  Symbolic Systems}.
\newblock {\em Inventiones Mathematicae}, 176(1):131--167, 2009.

\bibitem[Kit98]{kitchens1998sd}
B.~Kitchens.
\newblock {\em {Symbolic dynamics}}.
\newblock Springer New York, 1998.

\bibitem[LM95]{lind1995isd}
D.~Lind and B.~Marcus.
\newblock {\em {An Introduction to Symbolic Dynamics and Coding}}.
\newblock Cambridge University Press, 1995.

\bibitem[Moz89]{mozes1989tss}
S.~Mozes.
\newblock {Tilings, substitution systems and dynamical systems generated by
  them}.
\newblock {\em Journal d'analyse math{\'e}matique(Jerusalem)}, 53:139--186,
  1989.

\bibitem[PS10]{pavlovschraudner2009}
R.~Pavlov and M.~Schraudner.
\newblock Classification of sofic projective subdynamics of multidimensional
  shifts of finite type.
\newblock {\em Submitted}, 2010.

\bibitem[RJ87]{rogersjr1987trf}
H.~Rogers~Jr.
\newblock {\em {Theory of recursive functions and effective computability}}.
\newblock MIT Press Cambridge, MA, USA, 1987.

\bibitem[Rob71]{robinson1971uan}
R.M. Robinson.
\newblock {Undecidability and nonperiodicity for tilings of the plane}.
\newblock {\em Inventiones Mathematicae}, 12(3):177--209, 1971.

\end{thebibliography}

\end{document}